\documentclass[reqno]{amsart}

\usepackage{amsmath,amssymb,color,enumerate}

\newtheorem{theorem}{Theorem}[section]
\newtheorem{lemma}[theorem]{Lemma}
\newtheorem{proposition}[theorem]{Proposition}
\newtheorem{corollary}[theorem]{Corollary}

\numberwithin{equation}{section}

\newcommand{\mc}[1]{{\mathcal #1}}
\newcommand{\mf}[1]{{\mathfrak #1}}
\newcommand{\mb}[1]{{\mathbf #1}}
\newcommand{\bb}[1]{{\mathbb #1}}

\newcommand{\<}{\langle}
\renewcommand{\>}{\rangle}

\begin{document}

\title[Exclusion processes with conductances]{Hydrodynamic limit
  of gradient exclusion processes with conductances}

\author[T. Franco]{Tertuliano Franco}
\address{Tertuliano Franco \hfill\break\indent 
IMPA \hfill\break\indent 
Estrada Dona Castorina 110, \hfill\break\indent
J. Botanico, 22460 Rio
de Janeiro, Brazil} 
\email{tertu@impa.br}

\author[C. Landim]{Claudio Landim}
\address{Claudio Landim \hfill\break\indent 
IMPA \hfill\break\indent 
Estrada Dona Castorina 110, \hfill\break\indent
J. Botanico, 22460 Rio
de Janeiro, Brazil\hfill\break\indent 
{\normalfont and} \hfill\break\indent 
CNRS UMR 6085, \hfill\break\indent 
Avenue de l'Universit\'e, BP.12, Technop\^ole du Madrillet,
\hfill\break\indent  
F76801 Saint-\'Etienne-du-Rouvray, France.
} 
\email{landim@impa.br}


\noindent 


\begin{abstract}
  Fix a strictly increasing right continuous with left limits function
  $W: \bb R \to \bb R$ and a smooth function $\Phi : [l,r] \to \bb R$,
  defined on some interval $[l,r]$ of $\bb R$, such that $0<b \le
  \Phi'\le b^{-1}$. We prove that the evolution, on the diffusive
  scale, of the empirical density of exclusion processes, with
  conductances given by $W$, is described by the weak solutions of the
  non-linear differential equation $\partial_t \rho = (d/dx)(d/dW)
  \Phi(\rho)$. We derive some properties of the operator
  $(d/dx)(d/dW)$ and prove uniqueness of weak solutions of the
  previous non-linear differential equation.
\end{abstract}

\maketitle

\section{Introduction}
\label{sec1}

Recently attention has been raised to the hydrodynamic behavior of
interacting particle systems with random conductances \cite{n, jl2,
  fag, fjl}. In \cite{fjl}, for instance, the authors considered, for
a double sided $\alpha$--stable subordinator $W$, $0<\alpha<1$, the
nearest-neighbor one-dimensional exclusion process on $N^{-1} \bb Z$
in which a particle jumps from $x/N$ (resp. $(x+1)/N$) to $(x+1)/N$
(resp. $x/N$) at rate $\{N [W(x+1/N) - W(x/N)]\}^{-1}$. Their main
result can be restated as follows. On the the diffusive scale, as the
parameter $N\uparrow\infty$, the empirical density evolves according
to the solution of the differential equation
\begin{equation}
\label{aa1}
\partial_t \rho \;=\; \frac{d}{dx} \frac d{dW} \rho\;.
\end{equation}

The interesting feature is that, in constrast with usual
homogeneization phenomena, the entire noise survives in the limit and
the differential operator itself depends on the specific realization
of the Levy process $W$. The second surprising aspect is the
differential equation in which appears the derivative with respect to
a strictly increasing function $W$ which may have jumps. In fact, in
the Levy case, the set of discontinuities is dense on $\bb R$.

While the operators $(d/dW) (d/dx)$ have attracted much attention,
being closely related to the so-called gap diffusions or
quasi-diffusions when $W$ has no jumps \cite{m}, the operator $(d/dx)
(d/dW)$ have not been examined yet in the case where $W$ exhibit
jumps.  We refer to \cite{lo1, lo2, f} for recent results on the
operators $(d/dx) (d/dW)$ in the case where $W$ are increasing
continuous funcions.

As we shall see below, non-linear versions of the partial differential
equation \eqref{aa1} appear naturally as scaling limits of interacting
particle systems in inhomogeneous media. They may model diffusions in
which permeable membranes, at the points of the discontinuities of
$W$, tend to reflect particles, creating space discontinuities in the
solutions.

We present in this paper a gradient exclusion process whose
macroscopic evolution is described by the nonlinear differential
equation
\begin{equation*}
\partial_t \rho \;=\; \frac{d}{dx} \frac d{dW} \Phi(\rho)\;,
\end{equation*}
where $\Phi$ is a smooth function strictly increasing in the range of
$\rho$.  To prove this result, we examine in details the operator
$(d/dx) (d/dW)$ in $L^2(\bb T)$, where $\bb T$ is the one-dimensional
torus. We prove, in Theorem \ref{s11}, that $(d/dx) (d/dW)$, defined
on an appropriate domain, is non-positive, self-adjoint and
dissipative, and thus, the infinitesimal generator of a reversible
Markov process. We also prove that the eigenvalues of $-(d/dx) (d/dW)$
are countable and have finite multiplicity, the associated
eigenvectors forming a complete orhtonormal system.

\section{Notation and Results}
\label{sec2}

We examine the hydrodynamic behavior of a one-dimensional exclusion
process with conductances given by a strictly increasing function.
Let $\bb T_N$ be the one-dimensional dicrete torus with $N$
points. Distribute particles on $\bb T_N$ in such a way that each site
of $\bb T_N$ is occupied by at most one particle. Denote by $\eta$ the
configurations of state space $\{0,1\}^{\bb T_N}$ so that $\eta(x) =0$
if site $x$ is vacant and $\eta(x)=1$ if site $x$ is occupied.

Fix $a> -1/2$ and a \emph{strictly increasing} right continuous with
left limits (c\`adl\`ag) function $W: \bb R \to \bb R$, periodic in
the sense that $W(u+1) - W(u) = W(1) - W(0)$ for all $u$ in $\bb
R$. To simplify notation assume that $W$ vanishes at the origin,
$W(0)=0$. For $0\le x\le N-1$, let
\begin{equation*}
c_{x,x+1}(\eta) \;=\; 1 \;+\; a \{ \eta(x-1) + \eta(x+2)\}\;,
\end{equation*}
where all sums are modulo $N$, and let
\begin{equation*}
\xi_x \;=\; \frac 1{N[W(x+1/N) - W(x/N)]}
\end{equation*}
with the convention that $\xi_{N-1} = \{N[W(1) - W(1-[1/N])]\}^{-1}$.

The stochastic evolution can be described as follows.  At rate $\xi_x
c_{x,x+1}(\eta)$ the occupation variables $\eta(x)$, $\eta(x+1)$ are
exchanged. Note that if $W$ is differentiable at $x/N$, the rate at
which particles are exchanged is of order $1$, while if $W$ is
discontinuous, the rate is of order $1/N$. Assume, to fix ideas, that
$W$ is discontinuous at some point $x/N$ and smooth on intervals
$(x/N, x/N + \epsilon)$, $(x/N - \epsilon, x/N )$. In this case, the
rate at which particles jump over the bond $\{x-1, x\}$ is of order
$1/N$, while in a neighborhood of size $N$ of this bond, particles
jump at rate $1$.  Particles are thus reflected over the bond
$\{x-1,x\}$. However, since time will be scaled diffusively and since
on a time interval of length $N^2$ a particle spends a time of order
$N$ at site $x$, particle will be able to jump over the slower bond
$\{x-1, x\}$. This bond may model a membrane which difficults the
passage of particles.

The effect of the factor $c_{x,x+1}(\eta)$ is less dramatic. If the
parameter $a$ is positive, the presence of particles at the neighbor
sites of the bond $\{x,x+1\}$ speeds up the exchange by a factor of
order one.

The dynamics informally presented describes a Markov evolution. The
generator $L_N$ of this Markov process acts on functions $f:
\{0,1\}^{\bb T_N} \to \bb R$ as
\begin{equation}
\label{g4}
(L_N f) (\eta) \;=\; \sum_{x \in \bb Z} \xi_x\, c_{x,x+1}(\eta)\,
\{ f(\sigma^{x,x+1} \eta) - f(\eta) \} \;,
\end{equation}
where $\sigma^{x,x+1} \eta$ is the configuration obtained from $\eta$
by exchanging the variables $\eta(x)$, $\eta(x+1)$:
\begin{equation}
\label{g5}
(\sigma^{x,x+1} \eta)(y) \;=\;
\begin{cases}
\eta (x+1) & \text{ if } y=x,\\
\eta (x) & \text{ if } y=x+1,\\
\eta (y) & \text{ otherwise}.
\end{cases}
\end{equation}

A simple computation shows that the Bernoulli product measures
$\{\nu^N_\alpha : 0\le \alpha \le 1\}$ are invariant, in fact
reversible, for the dynamics. The measure $\nu^N_\alpha$ is obtained
by placing a particle at each site, independently from the other
sites, with probability $\alpha$. Thus, $\nu^N_\alpha$ is a product
measure over $\{0,1\}^{\bb T_N}$ with marginals given by
\begin{equation*}
\nu^N_\alpha \{\eta : \eta(x) =1\} \;=\; \alpha
\end{equation*}
for $x$ in $\bb T_N$. We wil often omit the index $N$ of
$\nu^N_\alpha$. 

Denote by $\{\eta_t : t\ge 0\}$ the Markov process on $\{0,1\}^{\bb
  T_N}$ associated to the generator $L_N$ \emph{speeded up} by
$N^2$. Let $D(\bb R_+, \{0,1\}^{\bb T_N})$ be the path space of
c\`adl\`ag trajectories with values in $\{0,1\}^{\bb T_N}$. For a
measure $\mu_N$ on $\{0,1\}^{\bb T_N}$, denote by $\bb P_{\mu_N}$ the
probability measure on $D(\bb R_+, \{0,1\}^{\bb T_N})$ induced by the
initial state $\mu_N$ and the Markov process $\{\eta_t : t\ge 0\}$.
Expectation with respect to $\bb P_{\mu_N}$ is denoted by $\bb
E_{\mu_N}$.

Denote by $\bb T$ the one-dimensional torus $[0,1)$.  A sequence of
probability measures $\{\mu_N : N\geq 1 \}$ on $\{0,1\}^{\bb T_N}$ is
said to be associated to a profile $\rho_0 :\bb T \to [0,1]$ if
\begin{equation}
\label{f09}
\lim_{N\to\infty}
\mu_N \Big\{ \, \Big\vert \frac 1N \sum_{x\in\bb T_N} H(x/N) \eta(x)
- \int H(u) \rho_0(u) du \Big\vert > \delta \Big\} \;=\; 0
\end{equation}
for every $\delta>0$ and every continuous functions $H: \bb T \to \bb
R$.

\subsection{The operator $\mc L_W$}

Denote by $\<\cdot , \cdot\>$ the inner product of $L^2(\bb T)$:
\begin{equation*}
\< f, g\>\;=\; \int_{\bb T} f(u)\, g(u)\, du\;.
\end{equation*}
Let $\mc D_W$ be the set of functions $f$ in $L^2(\bb T)$ such that
\begin{equation*}
f(x) \;=\; a \;+\; b W(x) \;+\; \int_{(0,x]} W(dy) \int_0^y \mf f(z) \, dz
\end{equation*}
for some function $\mf f$ in $L^2(\bb T)$ such that 
\begin{equation*}
\int_0^1 \mf f(z) \, dz \;=\; 0\;, \quad
\int_{(0,1]} W(dy) \Big\{ b + \int_0^y \mf f(z) \, dz \Big\} \;=\;0\; .
\end{equation*}
Define the operator $\mc L_W : \mc D_W \to L^2(\bb T)$ such that $\mc
L_W f = \mf f$. Denote by $\bb I$ the identity operator in $L^2(\bb
T)$. 

\begin{theorem}
\label{s11}
The operator $\mc L_W : \mc D_W \to \bb L^2(\bb T)$ enjoys the
following properties.

\renewcommand{\theenumi}{\alph{enumi}}
\renewcommand{\labelenumi}{{\rm (\theenumi)}}

\begin{enumerate}
\item $\mc D_W$ is dense in $\bb L^2(\bb T)$;

\item The operator $\bb I - \mc L_W : \mc D_W \to L^2(\bb T)$ is
  bijective. 

\item $\mc L_W: \mc D_W \to \bb L^2(\bb T)$ is self-adjoint and
  non-positive:
\begin{eqnarray*}
\< -\mc L_W f , f\> \;\ge\; 0
\end{eqnarray*}

\item $\mc L_W$ is dissipative;

\item The eigenvalues of the operator $- \mc L_W$ form a countable set
  $\{\lambda_n : n\ge 0\}$. All eigenvalues have finite multiplicity,
  $0= \lambda_0 \le \lambda_1 \le \cdots$, and $\lim_{n\to\infty} \lambda_n
  = \infty$.

\item The eigenvectors $\{f_n\}$ form a complete orhtonormal system
  
\end{enumerate}
\end{theorem}

In view of (a), (b), (d), by the Hille-Yosida theorem, $\mc L_W$ is
the generator of a strongly continuous contraction semi-group
semigroup $\{P_t : t\ge 0\}$, $P_t :L^2(\bb T) \to L^2(\bb T)$. Denote
by $\{G_\lambda : \lambda >0\}$, $G_\lambda :L^2(\bb T) \to L^2(\bb
T)$, the semi-group of resolvants associated to the operator $\mc
L_W$: $G_\lambda = (\lambda - \mc L_W)^{-1}$. In terms of the
semi-group $\{P_t\}$, $G_\lambda = \int_0^\infty e^{-\lambda t} P_t \,
dt$.

\subsection{The hydrodynamic equation}
\label{ss2.3}

For a positive integer $m\ge 1$, denote by $C^m(\bb T)$ the space of
continuous functions $H:\bb T\to \bb R$ with $m$ continuous
derivatives. Fix $l<r$ and a smooth function $\Phi :[l,r]\to \bb R$
whose derivative is bounded below by a strictly positive constant and
bounded above by a finite constant:
\begin{equation*}
0 \;<\, B^{-1} \le \Phi'(u)\; \le\; B
\end{equation*}
for $u$ in $[l,r]$. Consider a bounded density profile $\gamma : \bb T
\to [l,r]$. A bounded function $\rho : \bb R_+ \times \bb T \to [l,r]$
is said to be a weak solution of the parabolic differential equation
\begin{equation}
\label{g03}
\left\{
\begin{array}{l}
{\displaystyle \partial_t \rho \; =\; \mc L_W \Phi(\rho) } \\
{\displaystyle \rho(0,\cdot) \;=\; \gamma(\cdot)}
\end{array}
\right.
\end{equation}
if for all functions $H$ in $C^1(\bb T)$, all $t>0$ and all
$\lambda>0$, 
\begin{equation*}
\< \rho_t, G_\lambda H\> \;-\; \< \gamma , G_\lambda H\> 
\;=\; \int_0^t \< \Phi(\rho_s) , \mc L_W G_\lambda H \>\, ds\;.
\end{equation*}
We prove in Section \ref{sec6} uniqueness of weak solutions. Existence
follows from the tightness of the sequence of probability measures 
$\bb Q_{\mu_N}^{W,N}$ introduced in Section \ref{sec3}.

\begin{theorem}
\label{t01}
Fix a continuous initial profile $\rho_0 : \bb T \to [0,1]$ and
consider a sequence of probability measures $\mu_N$ on $\{0,1\}^{\bb
  T_N}$ associated to $\rho_0$. Then, for any $t\ge 0$,
\begin{equation*}
\lim_{N\to\infty}
\bb P_{\mu_N} \Big\{ \, \Big\vert \frac 1N \sum_{x\in\bb T_N} 
H(x/N) \eta_t(x) - \int H(u) \rho(t,u) du \Big\vert 
> \delta \Big\} \;=\; 0
\end{equation*}
for every $\delta>0$ and every continuous functions $H$. Here, $\rho$
is the unique weak solution of the non-linear equation \eqref{g03}
with $l=0$, $r=1$, $\gamma = \rho_0$ and $\Phi(\alpha) = \alpha + a
\alpha^2$.
\end{theorem}

Denote by $\pi^{N}_{t}$ the empirical measure at time $t$.  This is
the measure on $\bb T$ obtained by rescaling space by $N$ and by
assigning mass $N^{-1}$ to each particle:
\begin{equation*}
\pi^{N}_{t} \;=\; \frac{1}{N} \sum _{x\in \bb T_N} \eta_t (x)\,
\delta_{x/N}\;,
\end{equation*}
where $\delta_u$ is the Dirac measure concentrated on $u$.

Theorem \ref{t01} states that the empirical measure $\pi^N_t$
converges, as $N\uparrow\infty$, to an absolutely continuous measure
$\pi(t,du) = \rho(t,u) du$, whose density $\rho$ is the solution of
\eqref{g03}. In Sections \ref{sec3}, \ref{sec4} we prove that $\rho$
has finite energy: for all $t>0$,
\begin{equation*}
\int_0^t ds \int_{\bb T} \Big\{ \frac d{dW} \Phi(\rho(s,u)) \Big\}^2
\, dW \;<\; \infty\;.
\end{equation*}
The derivative $d/{dW} \Phi(\rho(s,u))$ must be understood in the
generalized sense. Details are given in Section \ref{sec4}.

\section{The operator $\mc L_W$}
\label{sec5}

We examine in this section properties of the operator $\mc L_W$
introduced in the previous section. Recall that we denote by $\<\cdot,
\cdot\>$ the inner product of the Hilbert space $L^2(\bb T)$ and by
$\Vert \cdot \Vert$ its norm.

Let $\bb D (f)$ be the set of discontinuity points of a function
$f:\bb T\to \bb R$. Denote by $C_{W}(\bb T)$ the set of c\`adl\`ag
functions $f:\bb T \to\bb R$ such that $\bb D(f) \subset \bb D (W)$.
$C_{W}(\bb T)$ is provided with the usual sup norm $\|\cdot
\|_\infty$.

All functions in $C_W(\bb T)$ are bounded. In fact, it is easy to
prove that for each fixed $f$ in $C_W(\bb T)$ and $\epsilon >0$, there
exists $n\ge 1$ and $0\le z_1 <z_2 <\cdots <z_{n}<1$ such that
\begin{equation}
\label{f11}
|f(x) - f(y)| \;\le\; \epsilon \text{ for all $z_k \le x,y <z_{k+1}$,
  $1\le k\le n$,}
\end{equation}
where $z_{n+1} = z_1$.

Define the generalized derivative $\frac{ d}{dW}$ as follows
\begin{equation}
\label{malditesta}
\frac{d f}{dW} (x) = \lim_{\epsilon\rightarrow 0} \frac{f(x+\epsilon)
-f(x)}{W(x+\epsilon) -W(x)}\;,
\end{equation}
if the above limit exists and is finite. Denote by $\mf D_W$ the set
of functions $f$ in $C_{W} (\bb T)$ such that $\frac{df}{dW}(x)$ is
well defined and derivable, and
$\frac{d}{dx}\bigl(\frac{df}{dW}\bigr)$ belongs to $C_{W}(\bb T)$.
Define the operator $\mf L_W : \mf D_W \to C_{W}(\bb T)$ by
\begin{equation*}
\mf L_W f \;=\; \frac{d}{dx} \frac{d}{dW}f \;=\;
\frac{d}{dx}\left(\frac{df}{dW}\right)\;.
\end{equation*}

By \cite[Lemma 0.9 in Appendix]{dyn2}, given a right continuous
function $f$ and a continuous function $h$, 
\begin{equation*}
\frac{df}{dW} (x) \;= \;h(x)
\end{equation*}
for all $x$ in $\bb T$ if and only if
\begin{equation}
\label{f12}
f(b)- f(a) \;=\; \int_{(a,b]} h(y) dW (y) 
\end{equation}
for all $a<b$. Note that the function $h$ has integral equal to zero,
$\int_{\bb T} h \, dW =0$, because $f(1)=f(0)$.

It follows from this observation and the definition of the operator
$\mf L_W$ that $\mf D_W$ is the set of functions $f$ in $C_{W}(\bb T)$
such that
\begin{equation}
\label{f17}
f(x) \;=\; a \;+\; b W(x)\; +\; \int_{(0,x]}  dW(y) \int_0^y g(z)
\, dz
\end{equation}
for some function $g$ in $C_{W}(\bb R)$ and two reals $a$, $b$ such
that 
\begin{equation}
\label{f14}
b W(1) \;+\; \int_{\bb T} dW(y) \int_0^y g(z) \, dz \;=\;0\;, \quad
\int_{\bb T} g(z) \, dz \;=\;0 \;.
\end{equation}
The first requirement corresponds to the boundary condition
$f(1)=f(0)$ and the second one to the boundary condition $(df/dW) (1)
= (df/dW) (0)$. Equivalently, \eqref{f14} follow from the conditions
\begin{equation}
\label{f18}
\int_{\bb T} \frac {df}{dW}\,  dW \;=\; 0\;, 
\quad \int_{\bb T} \frac {d}{dx} \frac {df}{dW}\, dx \;=\; 0\;.
\end{equation}
One can check that the function $g$, as well as the constants $a$,
$b$, are unique.

\begin{lemma}
\label{s07}
The following statements hold.  

\renewcommand{\theenumi}{\alph{enumi}}
\renewcommand{\labelenumi}{{\rm (\theenumi)}}

\begin{enumerate}
\item The set $\mf D_W$ is dense in $L^2(\bb T)$.

\item The operator $\mf L_W : \mf D_W \to L^2(\bb T)$ is symmetric and
nonpositive. More precisely,
\begin{equation*}
\< \mf L_W f, g\> \;=\; - \int_{\bb T} \frac {df}{dW} \, \frac
{dg}{dW}\, dW
\end{equation*}
for all $f$, $g$ in $\mf D_W$.

\item $\mf L_W$ satisfies a Poincar\'e inequality: There exists a
  finite constant $C_0$ such that
\begin{equation*}
\Vert f\Vert^2 \;\le\; C_0 \< - \mf L_W f, f\> \;+\; \Big(\int_{\bb T}
f(x) \, dx\Big)^2 
\end{equation*}
for all functions $f$ in $\mf D_W$.

\item The Green's function $G$ of $\mf L_W$ is given by
\begin{equation*}
G(x,y) \;=\; 
\left\{
\begin{array}{ll}
{\displaystyle - \frac{[W(y) - W(0)]\, [W(1) - W(x)]}{W(1)-W(0)} }& 0\le
    y\le x\le 1\;, \\ \\
{\displaystyle - \frac{[W(1) - W(y)]\, [W(x) - W(0)]}{W(1)-W(0)} }& 0\le
    x\le y\le 1\;.
\end{array}
\right.
\end{equation*}
\end{enumerate}
\end{lemma}

\begin{proof}
Since the continuous functions are dense in $L^2(\bb T)$, to prove (a)
it is enough to show that for each continuous function $f:\bb T\to \bb
R$ and $\epsilon >0$, there exists $g$ in $\mf D_W$ such that $\Vert
f-g\Vert \le \epsilon$.

Fix therefore a continuous function $f:\bb T\to \bb R$ and $\epsilon
>0$. There exists $\delta >0$ such that $|f(y)-f(x)|\le \epsilon$ if
$|x-y|\le \delta$. Choose an integer $n \ge \delta^{-1}$ and consider
the function $g:\bb T \to \bb R$ defined by
\begin{equation*}
g(x) \;=\; \sum_{j=0}^{n-1}\frac {f([j+1]/n) - f(j/N)}{W([j+1]/n) -
  W(j/N)} \mb 1\{(j/n, (j+1)/n]\}(x)\; , 
\end{equation*}
where $\mb 1\{A\}$ stands for the indicator of the set $A$. Let $G:
\bb T \to \bb R$ be given by $G(x) = f(0) + \int_{(0,x]} g(y)
W(dy)$. By definition of $g$, $G(j/n) = f(j/n)$ for $0\le j <
n$. Thus, by our choice of $n$ and by definition of $G$, for $j/n \le
x \le (j+1)/n$,
\begin{equation*}
\big| G(x) - f(x) \big| \;\le\; \big| G(x) - G(j/n) \big| \;+\; 
\big| f(x) - f(j/n) \big| \;\le\; 2 \epsilon\;.
\end{equation*}
so that $\Vert G - f \Vert_\infty \le 2 \epsilon$ if $\Vert
\cdot\Vert_\infty$ stands for the sup norm. Note that
\begin{equation}
\label{f15}
\int_{(0,1]} g\, dW \;=\;0\;.
\end{equation}

It remains to show that the function $G$ may be approximated in
$L^2(\bb T)$ by functions in the domain $\mf D_W$. Note that we were
free to choose the set $\{0, 1/n, \dots, (n-1)/n\}$ as long as the
distance between two consecutive points is bounded by $\delta$. We may
therefore assume, without loss of generality, that $W$ is continuous at
these points. Denote by $\{H_k : k\ge 1\}$ a sequence of smooth
functions $H_k : \bb T \to \bb R$ absolutely bounded by $\Vert
g\Vert_\infty$ and such that $\lim_k H_k(x) = g(x)$ for $xn \not\in
\bb Z$. By the dominated convergence theorem,
\begin{equation}
\label{f16}
\lim_{k\to\infty} \int_{\bb T} \big\vert H_k(y) - g(y) \big|\,
dW(y)\;=\; 0\;.
\end{equation}
Let $\{F_k : k\ge 1\}$ be the sequence of functions $F_k:\bb T\to\bb
R$ defined by
\begin{eqnarray*}
\!\!\!\!\!\!\!\!\!\!\!\!\!\!\! &&
F_k(x)\;=\; f(0) \;+\; \int_{(0,x]} \Big \{ b_k + \int_0^y H'_k(z) \, 
dz \Big \} W(dy) \\
\!\!\!\!\!\!\!\!\!\!\!\!\!\!\! && \qquad
\;=\; f(0) \;+\; b_k \, W(x) \;+\;
\int_{(0,x]} W(dy) \int_0^y H'_k(z) \, dz \;,
\end{eqnarray*}
where $b_k = H_k(0) - W(1)^{-1} \int_{(0,1]} H_k(y) \, dW(y)$. By
\eqref{f15}, \eqref{f16}, $F_k$ converges in the uniform topology to
$G$. On the other hand, in view of \eqref{f17} and our choice of
$b_k$, $F_k$ belongs to $\mf D_W$ for each $k\ge 1$ because $H'_k$,
being continuous, belongs to $C_W(\bb T)$.  This concludes the proof
of (a).

To prove (b), fix two functions $f$, $g$ in $\mf D_W$ and let $F=
df/dW$. $F$ is differentiable with derivative in $C_W(\bb T)$. Fix
$\epsilon >0$ and denote by $\{z_1, \dots , z_n\}$ the finite set
given by \eqref{f11} for the function $g$. Adding extra points if
necessary, we may assume that $\max_{1\le k\le n} \sup_{z_k\le x,y\le
  z_{k+1}} |F(y)-F(x)|\le \epsilon$ because $F$ is
continuous. Decomposing the integral over $\bb T$ on the intervals
$[z_k, z_{k+1}]$, we get that
\begin{equation*}
\< \mf L_W f, g\> \;=\; \int_{\bb T} \frac {dF}{dx} (x) \, g(x) \,
dx\; =\; \sum_{k=1}^n g(z_k) \{F(z_{k+1}) - F(z_k)\} \; +\; \epsilon
\, \Big\Vert  \frac {dF}{dx} \Big\Vert_\infty\;.
\end{equation*}
Changing the order of summations in the last term, in view of
\eqref{f12}, we obtain that the previous sum is equal to
\begin{equation*}
- \sum_{k=1}^n \{ g(z_k) - g(z_{k-1}) \} F(z_k) \;=\; 
- \sum_{k=1}^n F(z_k) \, \int_{(z_{k-1}, z_k]} \frac {dg}{dW}(x) \,
dW(x)\; .
\end{equation*}
Recall that $dg/dW$ is continuous and that $|F(x) -
F(z_k)|\le\epsilon$ for $z_{k-1} \le x\le z_k$. The previous sum is
thus equal to
\begin{equation*}
- \int_{\bb T} F(x)\, \frac {dg}{dW}(x) \, dW(x)\; +\; 
\epsilon \, \Big\Vert  \frac {dg}{dW} \Big\Vert_\infty\,
[W(1)-W(0)]\;. 
\end{equation*}
This proves the first identity from which it follows that $\mf L_W$ is
symmetric and nonpositive.

To prove the Poincar\'e inequality, fix a function $f$ in $\mf D_W$
and observe that by \eqref{f12}
\begin{eqnarray*}
\!\!\!\!\!\!\!\!\!\!\!\! &&
\int_{\bb T} f(x)^2 \, dx \;-\; \Big( \int_{\bb T} f(x)\, dx \Big)^2
\;=\; \int_{\bb T} \Big( \int_{\bb T} [f(x) - f(y)] \, dy \Big)^2 
\, dx \\
\!\!\!\!\!\!\!\!\!\!\!\! && \qquad\qquad
=\; \int_{\bb T} dx \Big( \int_{\bb T} dy \int_{(y,x]}
\frac {df}{dW} (z) \, dW(z) \Big)^2  \;.
\end{eqnarray*}
To conlude the proof, it remains to apply twice the Schwarz inequality
and to change the order of integration. Note that this proof gives
$C_0 = W(1) - W(0)$.

An elementary computation permits to check that the Green's function
is given by the expression proposed.
\end{proof}

Denote by $\<\cdot, \cdot\>_{1,2}$ the inner product on $\mf D_W$
defined by
\begin{equation*}
\<f,g\>_{1,2}\; =\; \<f,g\> \; +\; \<- \mf L_W f, g\>
\;=\; \<f,g\> \; +\; \int_{\bb T} \frac {df}{dW} \, \frac
{dg}{dW}\, dW\;.
\end{equation*}
Let $H^1_2 (\bb T)$ be the set of all functions $f$ in $L^2(\bb T)$
for which there exists a sequence $\{f_n : n\ge 1\}$ in $\mf D_W$
such that $f_n$ converges to $f$ in $L^2(\bb T)$ and $f_n$ is Cauchy
for the inner product $\<\cdot, \cdot \>_{1,2}$. Such sequence
$\{f_n\}$ is called admissible for $f$. For $f$, $g$ in
$H^1_2 (\bb T)$, define
\begin{equation}
\label{f20}
\<f,g\>_{1,2}\; =\; \lim_{n\to\infty} \<f_n,g_n\>_{1,2}\;,
\end{equation}
where $\{f_n\}$, $\{g_n\}$ are admissible sequences for $f$, $g$,
respectively.  By \cite[Proposition 5.3.3]{z}, this limit exists and
does not depend on the admissible sequence chosen. Moreover, $H^1_2
(\bb T)$ endowed with the scalar product $\<\cdot, \cdot \>_{1,2}$
just defined is a real Hilbert space.

Denote by $L^2_W(\bb T)$ the Hilbert space generated by the continuous
functions endowed with the inner product $\<\cdot, \cdot\>_W$ defined
by
\begin{equation*}
\<f,g\>_W \;=\; \int_{\bb T} f(x)\, g(x)\, W(dx)\;.
\end{equation*}
The norm associated to the scalar product $\<\cdot, \cdot \>_W$ is
denoted by $\Vert \cdot\Vert_W$.

\begin{lemma}
\label{s09}
A function $f$ in $L^2(\bb T)$ belongs to $H^1_2 (\bb T)$ if and only
if there exists $F$ in $L^2_W(\bb T)$ and a finite constant $c$ such that
\begin{equation*}
\int_{(0,1]} F(y) \, dW(y) \;=\;0 \quad\text{and}\quad
f(x) \;=\; c\;+\; \int_{(0,x]} F(y) \, dW(y)
\end{equation*}
Lebesgue almost surely. We denote the generalized $W$-derivative $F$
of $f$ by $df/dW$. For $f$, $g$ in $H^1_2 (\bb T)$,
\begin{equation*}
\<f,g\>_{1,2}\; =\; \<f,g\> \;+\; \int_{\bb T} \frac{df}{dW}\,
\frac{dg}{dW}\, dW\;.
\end{equation*}
\end{lemma}

\begin{proof}
Fix $f$ in $H^1_2 (\bb T)$. By definition, there exists a sequence
$\{f_n : n\ge 1\}$ in $\mf D_W$ which converges to $f$ in $L^2(\bb T)$
and which is Cauchy in $H^1_2 (\bb T)$. In particular, $df_n/dW$ is
Cauchy in $L^2_W(\bb T)$ and therefore converges to some function $G$
in $L^2_W(\bb T)$. By \eqref{f18}, 
\begin{equation*}
\int_{\bb T} \frac {df_n}{dW}\,  dW \;=\; 0
\end{equation*}
for all $n\ge 1$ so that $\int_{(0,1]} G\, dW =0$. Let $g(x) =
\int_{(0,x]} G(y) dW(y)$. Since $\mb 1\{(x,y]\}$ belongs to $L^2_W(\bb
T)$, for all $x$, $y$ in $\bb T$,
\begin{equation*}
g(y) - g(x) \;=\; \int_{(x,y]} G \, dW
\;=\; \lim_{n\to\infty} \int_{(x,y]} \frac{df_n}{dW}\, dW
\;=\; \lim_{n\to\infty} \{f_n(y) - f_n(x)\} \;.
\end{equation*}

We claim that $\int_{\bb T} \{f_n(y) - f_n(x)\} dx$ converges to
$\int_{\bb T} \{g(y) - g(x)\} dx$ for all $y$ in $\bb T$. Indeed, on
the one hand, for each fixed $y$, $f_n(y) - f_n(x)$ converges to $g(y)
- g(x)$. On the other hand, by Schwarz inequality, 
\begin{equation*}
[f_n(y) - f_n(x)]^2 \;\le\; [W(1)-W(0] \int_{\bb T}
\Big(\frac{df_n}{dW}\Big)^2 \, dW \;\le\; C_0
\end{equation*}
for some finite constant $C_0$. It remains to apply the dominated
convergence theorem to conclude.

Since $f_n$ converges to $f$ in $L^2(\bb T)$, $\int_{\bb T} f_n(x)\,
dx$ converges to $\int_{\bb T} f(x)\, dx$. By Schwarz inequality, $g$
belongs to $L^2(\bb T)$ so that $\int_{\bb T} g(x) dx$ is
finite. Therefore, for all $y$ in $\bb T$,
\begin{eqnarray*}
\!\!\!\!\!\!\!\!\!\!\!\!\!\!\! &&
\lim_{n\to\infty}  f_n(y) \;=\; \lim_{n\to\infty} \Big\{
f_n(y) - \int_{\bb T}  f_n(x)\, dx \Big\} \;+\; \int_{\bb T} f(x)\, dx
\\
\!\!\!\!\!\!\!\!\!\!\!\!\!\!\! && \qquad
\;=\; g(y) - \int_{\bb T} g(x)\, dx  \;+\; \int_{\bb T} f(x)\, dx\;.
\end{eqnarray*}
Thus $f_n$ converges pointwisely to the above function.  As $f_n$ also
converges to $f$ in $L^2(\bb T)$, we deduce that $f = c + g$ a.s., and
thus in $L^2(\bb T)$, for $c= \int_{\bb T} f(x)\, dx - \int_{\bb T}
g(x)\, dx$, proving the first statement of the lemma.

The reciproque is simpler. Let $f = c + \int_{(0,x]} F(y) \, dW(y)$
for some $F$ in $L^2_W(\bb T)$ such that $\int_{(0,1]} F(y)\,
dW(y)=0$. There exists a sequence $\{g_n : n\ge 1\}$ of smooth
functions converging to $F$ in $L^2_W(\bb T)$ and such that
$\int_{(0,1]} g_n(y)\, dW(y)=0$. Let $f_n (x) = c + \int_{(0,x]} dW(y)
\{ g_n(0) + \int_0^y g'_n(z) dz\}$. For each $n\ge 1$, $f_n$ belongs
to $\mf D_W$ because $g'_n$ is continuous.  Schwarz inequality shows
that $f_n$ converges to $f$ in $L^2(\bb T)$. Finally, $\{f_n : n\ge
1\}$ is a Cauchy sequence for the inner product $\<\cdot, \cdot
\>_{1,2}$ because $df_n/dW = g_n$ converges to $F$ in $L^2_W(\bb T)$.
Note that we just proved that the sequence $\{f_n : n\ge 1\}$ is
admissible for $f$

Fix $f$, $g$ in $H^1_2(\bb T)$ and recall that we denote by $df/dW$,
$dg/dW$ the generalized $W$-derivatives of $f$, $g$, respectively.
Denote by $\{f_n : n\ge 1\}$, $\{g_n : n\ge 1\}$ the admissible
sequences constructed in the previous paragraph for $f$ and $g$,
respectively. By definition,
\begin{equation*}
\<f,g\>_{1,2} \;=\; \lim_{n\to\infty} \<f_n,g_n\>_{1,2}
\;=\; \lim_{n\to\infty} \Big\{ \<f_n,g_n\> \; +\; \int_{\bb T}
\frac{df_n}{dW}\, \frac{dg_n}{dW}\, dW \Big\} \;.
\end{equation*}
Since $f_n$ (resp. $g_n$) converges to $f$ (resp. $g$) in $L^2(\bb T)$
and since $df_n/dW$ (resp. $dg_n/dW$) converges to $df/dW$
(resp. $dg/dW$) in $L^2_W(\bb T)$, the previous expression is equal to
\begin{equation*}
\<f,g\> \; +\; \int_{\bb T} \frac{df}{dW}\, \frac{dg}{dW}\, dW\;.
\end{equation*}
This concludes the proof of the lemma.
\end{proof}

\begin{lemma}
\label{s08}
The embedding $H^1_2 (\bb T) \subset L^2(\bb T)$ is compact.
\end{lemma}

\begin{proof}
Consider a sequence $\{u_n : n\ge 1\}$ bounded in $H^1_2 (\bb
T)$. We need to prove the existence of a subsequence $\{u_{n_k} : k\ge
1\}$ which converges in $L^2(\bb T)$.


By the previous lemma, $u_n(x) = c_n + \int_{(0,x]} U_n (y)\, dW(y)$
for some $U_n$ in $L^2_W(\bb T)$ such that $\int_{(0,1]} U_n (y)\,
dW(y) =0$. Moreover, $\Vert U_n \Vert_W \le \Vert u_n\Vert_{1,2}$. The
sequence $\{U_n\}$ is therefore bounded in $L^2_W(\bb T)$. Also, by
Schwarz inequality, the sequence $\int_{(0,x]} U_n (y) dW(y)$ is
bounded in $L^2(\bb T)$. Since $c_n = u_n(x) - \int_{(0,x]} U_n (y)
dW(y)$ and since both sequence of functions on the right hand side are
bounded in $L^2(\bb T)$, the sequence of real numbers $\{c_n\}$ is
also bounded.

Since $\{U_n\}$ is a bounded sequence in $L^2_W(\bb T)$ and since the
sequence of real numbers $\{c_n\}$ is bounded, there exists a
subsequence $\{n_k\}$ such that $c_{n_k}$ converges and $U_{n_k}$
converges weakly in $L^2_W(\bb T)$ to a limit denoted by $U$. As
constants belong to $L^2_W(\bb T)$, $\int_{(0,1]} U (y) \, dW(y) =
\lim_k \int_{(0,1]} U_{n_k} (y) \, dW(y) =0$.  Moreover, for all $x$
in $\bb T$, as $\mb 1\{(0,x]\}$ belongs to $L^2_W(\bb T)$,
\begin{equation*}
\lim_{k\to\infty} u_{n_k}(x) \;=\; \lim_{k\to\infty} \Big\{ c_{n_k}
\;+\; \int_{(0,x]} U_{n_k} (y) dW(y) \Big\}
\;=\;  c \;+\; \int_{(0,x]} U (y) dW(y)\;, 
\end{equation*}
if $c$ stands for the limit of the sequence $c_{n_k}$. The sequence
$u_{n_k}$ thus converges pointwisely to $u(x) = c+ \int_{(0,x]} U (y)
dW(y)$. Since, by Schwarz inequality, $u_{n_k} (x)^2$ is bounded by
$2c_{n_k}^2 + 2 [W(1) - W(0)] \, \Vert U_{n_k}\Vert^2_W$, by the
dominated convergence theorem, $u_{n_k}$ converges to $u$ in $L^2(\bb
T)$. Note that the limit $u$ belongs to $H^1_2(\bb T)$.
\end{proof}

Let $\mc D_W$ be the set of functions $f$ in $H^1_2(\bb T)$ for which
there exists $u$ in $L^2(\bb T)$ such that
\begin{equation}
\label{f19}
\<f,g\>_{1,2}\; =\;
\<f,g\> \; + \; \int \frac {df}{dW}\, \frac {dg}{dW}\, dW \;=\; 
\<u,g\>
\end{equation}
for all $g$ in $H^1_2(\bb T)$. By Lemma \ref{s07} (b), $\mf D_W
\subset \mc D_W$ and, by definition, $\mc D_W \subset H^1_2(\bb
T)$. The function $u$ is uniquely determined because, by Lemma
\ref{s07} (a), $H^1_2(\bb T) \supset \mf D_W$ is dense in $L^2(\bb
T)$. By definition of $H^1_2(\bb T)$ and by \eqref{f20}, it is enough
to check \eqref{f19} for functions $g$ in $\mf D_W$.

\begin{lemma}
\label{s10}
The domain $\mc D_W$ consist of all functions $f$ in $L^2(\bb T)$ such
that
\begin{equation*}
f(x) \;=\; a \;+\; b W(x) \;+\; \int_{(0,x]} W(dy) \int_0^y \mf f(z) \, dz
\end{equation*}
for some function $\mf f$ in $L^2(\bb T)$ such that 
\begin{equation*}
\int_0^1 \mf f(z) \, dz \;=\; 0\;, \quad
\int_{(0,1]} W(dy) \Big\{ b + \int_0^y \mf f(z) \, dz \Big\} \;=\;0\; .
\end{equation*}
Moreover, in this case,
\begin{eqnarray*}
-\, \int \frac {df}{dW}\, \frac {dg}{dW}\, dW \;=\; 
\<\mf f ,g\>
\end{eqnarray*}
for all $g$ in $H^1_2(\bb T)$.
\end{lemma}

\begin{proof}
We first show that any function $f$ in $L^2(\bb T)$ with the
properties listed in the statement of the lemma belongs to $\mc
D_W$. Fix such a function and consider a sequence $\{\mf f_n : n\ge
1\}$ of smooth functions $\mf f_n : \bb T\to\bb R$ which converges to
$\mf f$ in $L^2(\bb T)$ and such that $\int_0^1 \mf f_n(z) \, dz =
0$. Let
\begin{equation*}
f_n(x) \;=\; a \;+\; \int_{(0,x]} W(dy) \Big\{ b_n \;+\; 
\int_0^y \mf f_n(z) \, dz \Big\} \; ,
\end{equation*}
where $b_n$ is chosen so that $\int_{(0,1]} W(dy) \{ b_n + \int_0^y
\mf f_n(z) \, dz\} =0$. Note that $f_n$ belongs to $\mf D_W$ for each
$n\ge 1$.

As $n\uparrow\infty$, $b_n$ converges to $b$, $f_n$ converges to $f$
in $L^2(\bb T)$ and $\{f_n\}$ is Cauchy for the $\Vert
\cdot\Vert_{1,2}$ norm. Thus, $f$ belongs to $H^1_2(\bb T)$ and
$\{f_n\}$ is an admissible sequence for $f$.

Fix $g$ in $\mf D_W$. We claim that
\begin{equation*}
\<f,g\>_{1,2} \;=\; \<f,g\> \;-\; \<\mf f \, ,g\>\;.
\end{equation*}
Indeed, as $g$ belongs to $\mf D_W$, by \eqref{f20}, $\<f,g\>_{1,2} =
\lim_n \<f_n,g\>_{1,2}$ because the sequence $\{g_n : n\ge 1\}$
constant equal to $g$ is admissible for $g$. By definition of the
inner product $\<\cdot, \cdot\>_{1,2}$ and since $\mf L_W f_n = \mf
f_n$, $\<f_n,g\>_{1,2} = \<f_n, g\> + \<- \mf L_W f_n, g\> = \<f_n,
g\> + \<- \mf f_n, g\>$. Since $f_n$, $\mf f_n$ converge in $L^2(\bb
T)$ to $f$, $\mf f$, respectively, the claim is proved. In particular
\eqref{f19} holds with $u=f-\mf f$. This proves that $f$ belongs to
$\mc D_W$ and the identity claimed.

Conversely, assume that $f$ belongs to $\mc D_W$ and satisfy
\eqref{f19} for some $u$ in $L^2(\bb T)$. Thus, there exists $v$
(equal to $f-u$) in $L^2(\bb T)$ such that
\begin{equation}
\label{f21}
- \int \frac {df}{dW}\, \frac {dg}{dW}\, dW \;=\; \<v,g\>
\end{equation}
for all $g$ in $\mf D_W$. Taking $g=1$ in this equation we obtain
that $\int_0^1 v(x) \, dx=0$.

Since $f$ belongs to $H^1_2(\bb T)$, by Lemma \ref{s09}, $f(x) = c +
\int_{(0,x]} F(y) dW(y)$ for some function $F$ in $L^2_W(\bb T)$ such
that $\int_{(0,1]} F(y) dW(y) =0$. To prove the lemma we need to show
that
\begin{equation*}
F(y) \; =\; b \;+\; \int_0^y \mf f(z) \, dz
\end{equation*}
for some finite constant $b$ and some function $\mf f$ in $L^2(\bb T)$
such that $\int_0^1 \mf f(z) \, dz = 0$.

Fix $g$ in $\mf D_W$ so that
\begin{equation*}
g(x) \;=\; a \;+\; \int_{(0,x]} G(y)\, dW(y) 
\end{equation*}
for some continuous function $G: \bb T \to \bb R$ such that $\int_0^1
G(y) \, dW(y) =0$.  Since the integral of $v$ (resp. $G$) with respect
to the Lebesgue measure (resp. the measure $dW$) vanishes, changing
the order of integration, we obtain that
\begin{equation*}
\int_0^1 v(x) \, g(x) \, dx\; =\; - \int_{(0,1]}  G(y) \, 
\int_0^y v(x) \, dx \; dW(y) \;.
\end{equation*}
Therefore, in view of \eqref{f21},
\begin{equation*}
\int_{(0,1]} G(y) \, \int_0^y v(x) \, dx \; dW(y) \;=\;
\int_{(0,1]} G(y) \, F(y) \,  dW(y) 
\end{equation*}
for all functions $g$ in $\mf D_W$.  The proof of Lemma \ref{s07} (a)
shows that the set $\{dg/dW : g\in \mf D_W\}$ is dense in $L^2_{W,0} =
\{H \in L^2_W(\bb T) : \int H dW =0\}$. In particular, $F(y) = c +
\int_0^y v(x) \, dx$ for some finite constant $c$. This concludes the
proof of the lemma.
\end{proof}

Recall that we denote by $\bb I$ the identity in $L^2(\bb T)$. By
Lemma \ref{s07}, the symmetric operator $(\bb I - \mf L_W) : \mf D_W
\to \bb L^2(\bb T)$, is strongly monotone:
\begin{eqnarray*}
  \< (\bb I - \mf L_W) f, f\> \;\ge\; \<f,f\>
\end{eqnarray*}
for all $f$ in $\mf D_W$. Denote by $\mc A_1 : \mc D_W\to \bb L^2(\bb
T)$ its Friedrichs extension, defined as $\mc A_1 f =u$, where $u$ is
the function in $L^2(\bb T)$ given by \eqref{f19}.  By \cite[Theorem
5.5.a]{z}, $\mc A_1$ is self-adjoint, bijective and
\begin{equation}
\label{f22}
\< \mc A_1 f, f\> \;\ge\; \<f,f\>
\end{equation}
for all $f$ in $\mc D_W$. Note that the Friedrichs extension of the
strongly monotone operator $(\lambda \bb I - \mf L_W)$, $\lambda>0$,
is $\mc A_\lambda = (\lambda - 1)\bb I + \mc A_1:\mc D_W \to \bb
L^2(\bb T)$.

Define $\mc L_W : \mc D_W \to L^2(\bb T)$ by $\mc L_W = \bb I - \mc
A_1$. In view of \eqref{f19}, $\mc L_Wf =u$ if and only if
\begin{eqnarray*}
-\, \int \frac {df}{dW}\, \frac {dg}{dW}\, dW \;=\;  \<u,g\>
\end{eqnarray*}
for all $g$ in $H^1_2(\bb T)$. In particular by Lemma \ref{s07} (b)
$\mc L_W f = \mf L_W f$ for all $f$ in $\mf D_W$. Moreover, if a
function $f$ in $\mc D_W$ is represented as in Lemma \ref{s08}, $\mc
L_W f = \mf f$. This identity together with the identification of the
space $\mc D_W$ provides the alternative definition of the operator
$\mc L_W$ presented just before the statement of Theorem \ref{s11}.

\begin{proof}[Proof of Theorem \ref{s11}]
It follows from Lemma \ref{s07} (a) that the domain $\mc D_W$ is dense
in $L^2(\bb T)$ because $\mf D_W \subset \mc D_W$. This proves (a).

By definition, $\bb I-\mc L_W = \mc A_1 : \mc D_W \to
\bb L^2(\bb T)$, which have been shown to be bijective. This proves
(b). 

The self-adjointness of $\mc L_W: \mc D_W \to L^2(\bb T)$ follows from
the one of $\mc A_1$ and the definition of $\mc L_W$ as $\bb I - \mc
A_1$. Moreover, from \eqref{f22} we obtain that $\< -\mc L_W f , f\>
\ge 0$ for all $f$ in $\mc D_W$.

To prove (d), fix a function $g$ in $\mc D_W$, $\lambda >0$ and let $f
= (\lambda \bb I - \mc L_W) g$. Taking the scalar product with respect
to $g$ on both sides of this equation, we obtain that
\begin{eqnarray*}
\lambda \< g , g\> \;+\; \< -\mc L_W g , g\>
\;=\; \< g , f\> \;\le\;  \< g , g\>^{1/2} \,  
\<f , f\>^{1/2}\;.
\end{eqnarray*}
Since $g$ belongs to $\mc D_W$, by (c), the second term on the left
hand side is positive. Thus, $\Vert \lambda g\Vert \le \Vert f \Vert =
\Vert (\lambda \bb I - \mc L_W) g \Vert $.

We have already seen that the operator $(\bb I - \mf L_W): \mf D_W \to
L^2(\bb T)$ is symmetric and strongly monotone.  By Lemma \ref{s08},
the embedding $H^1_2(\bb T) \subset L^2(\bb T)$ is compact. Therefore,
by \cite[Theorem 5.5.c]{z}, the Friedrichs extension of $(\bb I - \mf
L_W)$, denoted by $\mc A_1 : \mc D_W \to L^2(\bb T)$, satisfies claims
(e) and (f) with $1\le \lambda_1 \le \lambda_2 \le \cdots$,
$\lambda_n\uparrow\infty$. In particular, the operator $- \mc L_W =
\mc A_1 - \bb I$ has the same property with $0\le \lambda_1 \le
\lambda_2 \le \cdots$, $\lambda_n\uparrow\infty$. Since $0$ is an
eigenvalue of $-\mc L_W$ associated at least to the constants, (e) and
(f) are in force.
\end{proof}

It follows also from \cite[Theorem 5.5.c]{z} that $f_n$ belongs to
$H^1_2(\bb T)$ for all $n$.

\subsection{Random walk with conductances}

Recall that $\bb T_N$ stands for the discrete one-dimensional torus
with $N$ points and recall the definition of the sequence $\{\xi_x :
0\le x\le N-1\}$. Consider the random walk $\{X^N_t : t\ge 0\}$ on
$N^{-1} \bb T_N$ which jumps over the bond $\{x/N, (x+1)/N\}$ at rate
$N^2 \xi_x = N/\{W(x+1/N) - W(x/N)\}$. The generator $\bb L_N$ of
this Markov process writes
\begin{equation*}
(\bb L_N f)(x/N) \; =\; N^2 \xi_x \{ f(x+1/N) - f(x/N)\} \; +\; 
N^2 \xi_{x-1} \{f(x-1/N) - f(x/N)\} \;.
\end{equation*}

The counting measure $m_N$ on $N^{-1} \bb T_N$ is reversible for this
process. Denote by $\{P^N_t : t\ge 0\}$ (resp. $\{G_\lambda^N :
\lambda >0\}$) the semigroup (resp. the resolvent) associated to the
generator $\bb L_N$:
\begin{equation*}
G_\lambda^N H  \;=\; \int_0^\infty dt\, e^{-\lambda t}  P_t^N H
\end{equation*}
for $H: N^{-1}\bb T_N \to \bb R$.

Fix a function $H: N^{-1}\bb T_N \to \bb R$. For $\lambda >0$, let
$H_\lambda^N = G_\lambda^N H$ be the solution of the resolvant
equation
\begin{equation*}
\lambda H_\lambda^N \;-\; \bb L_N H_\lambda^N \;=\; H\;.
\end{equation*}
Taking the scalar product on both sides of this equation with respect
to $H_\lambda^N$, we obtain that for all $N\ge 1$
\begin{equation}
\label{f05}
\begin{array}{l}
{\displaystyle 
\frac{1}{N}\sum_{x \in \bb T_N}  
H_\lambda^N (x/N)^2 \;\leq\; \frac{1}{\lambda^2}
\frac{1}{N}\sum_{x \in \bb T_N}  H (x/N)^2 \;,}\\
{\displaystyle 
\quad \frac{1}{N}\sum_{x \in \bb T_N} \xi_x (\nabla_N
H_\lambda^N)(x/N)^2 \;\leq\; \frac{1}{\lambda}
\frac{1}{N}\sum_{x \in \bb T_N}  H (x/N)^2 \;,}
\end{array}
\end{equation}
where $\nabla_N$ stands for the discrete derivative: $(\nabla_N
H)(x/N) = N[H(x+1/N) - H(x/N)]$.

On the other hand, if $H:\bb T\to \bb R$ is a continuous function and
we denote also by $H$ its restriction to $N^{-1} \bb T_N$, by
\cite[Lemma 4.6]{fjl},
\begin{equation}
\label{f23}
\lim_{\lambda\to\infty} \limsup_{N\to\infty}
\frac{1}{N} \sum_{x \in \bb T_N} \big|\lambda H_\lambda^N (x/N)
- H (x/N)\big|\;.
\end{equation}
Note that in \cite{fjl}, the function $W$ is of pure jump type, while
here it is any strictly increasing c\`adl\`ag function. One can check,
however, that the proof applies to our general
case.  

It follows from \cite[Lemma 4.5 (iii)]{fjl} that for every continuous
function $H:\bb T\to \bb R$,
\begin{equation}
\label{f24}
\limsup_{N\to\infty}
\frac{1}{N} \sum_{x \in \bb T_N} \big| (G_\lambda^N H) (x/N)
- ( G_\lambda H) (x/N)\big|\;.
\end{equation}

\section{Scaling limit}
\label{sec3}

Let $\mc M$ be the space of positive measures on $\bb T$ with total
mass bounded by one endowed with the weak topology. Recall that
$\pi^{N}_{t} \in \mc M$ stands for the empirical measure at time $t$.
This is the measure on $\bb T$ obtained by rescaling space by $N$ and
by assigning mass $N^{-1}$ to each particle:
\begin{equation}
\label{f01}
\pi^{N}_{t} \;=\; \frac{1}{N} \sum _{x\in \bb T_N} \eta_t (x)\,
\delta_{x/N}\;,
\end{equation}
where $\delta_u$ is the Dirac measure concentrated on $u$. For a
continuous function $H:\bb T \to \bb R$, $\<\pi^N_t, H\>$ stands for
the integral of $H$ with respect to $\pi^N_t$:
\begin{equation*}
\<\pi^N_t, H\> \;=\; \frac 1N \sum_{x\in\bb T_N}
H (x/N) \eta_t(x)\;.
\end{equation*}
This notation is not to be confounded with the inner product in
$L^2(\bb T)$ introduced earlier. Also, when $\pi_t$ has a density
$\rho$, $\pi(t,du) = \rho(t,u) du$, we sometimes write $\<\rho_t, H\>$
for $\<\pi_t, H\>$.

Fix $T>0$. Let $D([0,T], \mc M)$ be the space of $\mc M$-valued
c\`adl\`ag trajectories $\pi:[0,T]\to\mc M$ endowed with the
\emph{uniform} topology.  For each probability measure $\mu_N$ on
$\{0,1\}^{\bb T_N}$, denote by $\bb Q_{\mu_N}^{W,N}$ the measure on
the path space $D([0,T], \mc M)$ induced by the measure $\mu_N$ and
the process $\pi^N_t$ introduced in \eqref{f01}.

Fix a continuous profile $\rho_0 : \bb T \to [0,1]$ and consider a
sequence $\{\mu_N : N\ge 1\}$ of measures on $\{0,1\}^{\bb T_N}$
associated to $\rho_0$ in the sense \eqref{f09}. Let $\bb Q_{W}$ be
the probability measure on $D([0,T], \mc M)$ concentrated on the
deterministic path $\pi(t,du) = \rho (t,u)du$, where $\rho$ is the
unique weak solution of \eqref{g03} with $\gamma = \rho_0$, $l=0$,
$r=1$ and $\Phi(\alpha) = \alpha + a\alpha^2$.

\begin{proposition}
\label{s15}
As $N\uparrow\infty$, the sequence of probability measures $\bb
Q_{\mu_N}^{W,N}$ converges in the uniform topology to $\bb Q_{W}$.
\end{proposition}

The proof of this result is divided in two parts. In Subection
\ref{ss1}, we show that the sequence $\{\bb Q_{\mu_N}^{W,N} : N\ge
1\}$ is tight and in Subection \ref{ss2} we characterize the limit
points of this sequence. 

\begin{proof}[Proof of Theorem \ref{t01}]
Since $\bb Q_{\mu_N}^{W,N}$ converges in the uniform topology to $\bb
Q_{W}$, a measure which is concentrated on a deterministic path, for
each $0\le t\le T$ and each continuous function $H:\bb T\to \bb R$,
$\<\pi^N_t, H\>$ converges in probability to $\int_{\bb T} du \,
\rho(t,u)$ $H(u)$, where $\rho$ is the unique weak solution of
\eqref{g03} with $l=0$, $r=1$, $\gamma=\rho_0$ and $\Phi(\alpha) =
\alpha + a \alpha^2$.
\end{proof}

\subsection{Tightness}
\label{ss1}

Tightness of the sequence $\{\bb Q_{\mu_N}^{W,N} : N\ge 1\}$ is proved
as in \cite{jl2, fjl}. by considering first the auxiliary $\mc
M$-valued Markov process $\{\Pi^{\lambda,N}_t : t\ge 0\}$, $\lambda>0$,
defined by
\begin{equation*}
\Pi^{\lambda,N}_t (H) \;=\; \< \pi^N_t, G_\lambda^N H\>\;=\;
\frac{1}{N} \sum _{x\in \bb Z} \bigl (G_\lambda^N H \bigr)(x/N)
\eta_t (x)\;,
\end{equation*}
$H$ in $C(\bb T)$, where $\{G_\lambda^N : \lambda >0\}$ is the
resolvent associated to the random walk $\{X^N_t : t\ge 0\}$
introduced in Section \ref{sec5}.

We first prove tightness of the process $\{\Pi^{\lambda,N}_t : 0\le t
\le T\}$ for every $\lambda>0$ and then show that $\{\Pi^{\lambda,N}_t
: 0\le t \le T\}$ and $\{\pi^{N}_t : 0\le t \le T\}$ are not far
appart if $\lambda$ is large.

In constrast with \cite{jl2,fjl}, $D([0,T], \mc M)$ is here endowed
with the uniform topology. It is well known \cite{kl} that to prove
tightness of $\{\Pi^{\lambda,N}_t : 0\le t \le T\}$ it is enough to
show tightness of the real-valued processes $\{\Pi^{\lambda,N}_t (H) :
0\le t \le T\}$ for a set of smooth functions $H:\bb T\to \bb R$ dense
in $C(\bb T)$ for the uniform topology.

Fix a smooth function $H: \bb T \to \bb R$. Denote by the same symbol
the restriction of $H$ to $N^{-1} \bb T_N$. Let $H_\lambda^N =
G_\lambda^N H$ so that
\begin{equation}
\label{f04}
\lambda H_\lambda^N \;-\; \bb L_N H_\lambda^N \;=\; H\;.
\end{equation}
Keep in mind that $\Pi^{\lambda,N}_t (H) = \<\pi^N_t, H_\lambda^N \>$
and denote by $M^{N,\lambda}_t$ the martingale defined by
\begin{equation}
\label{f10}
M^{N,\lambda}_t \;=\;  \Pi^{\lambda,N}_t (H) \;-\; 
\Pi^{\lambda,N}_0 (H) \;-\; \int_0^t ds \, N^2 L_N \<\pi^N_s ,
H_\lambda^N \> \;.
\end{equation}
Clearly, tightness of $\Pi^{\lambda,N}_t (H)$ follows from tightness
of the martingale $M^{N,\lambda}_t$ and tightness of the additive
functional $\int_0^t ds \, N^2 L_N \<\pi^N_s , H_\lambda^N \>$.

An elementary computation shows that the quadratic variation
$\<M^{N,\lambda}\>_t$ of the martingale $M^{N,\lambda}_t$ is given by
\begin{equation*}
\frac 1{N^2} \sum_{x\in \bb T_N} \xi_x \, [(\nabla_N H_\lambda^N)(x/N)]^2
\int_0^t c_{x,x+1}(\eta_s) \, [\eta_s(x+1) - \eta_s(x)]^2 \, ds\;.
\end{equation*}
In particular, by \eqref{f05},
\begin{equation*}
  \<M^{N,\lambda}\>_t \;\le\; \frac{C_0 t}{N^2} \sum_{x\in \bb T_N}
  \xi_x \, [(\nabla_N H_\lambda^N)(x/N)]^2 \;\le\; \frac{C (H)
  t}{\lambda N}
\end{equation*}
for some finite constant $C(H)$ which depends only on $H$. Thus, by
Doob inequality, for every $\lambda>0$, $\delta>0$ 
\begin{equation}
\label{f02}
\lim_{N\to\infty} \bb P_{\mu_N} \Big[ \sup_{0\le t\le T}
\big\vert M^{N,\lambda}_t \big\vert \, > \, \delta \Big] 
\;=\; 0\;.
\end{equation}
In particular, the sequence of martingales $\{M^{N,\lambda}_t : N\ge
1\}$ is tight for the uniform topology.

It remains to examine the additive functional of the decomposition
\eqref{f10}. A long and elementary computations shows that $N^2 L_N
\<\pi^N , H_\lambda^N \>$ is equal to
\begin{eqnarray*}
\!\!\!\!\!\!\!\!\!\!\!\!\!\! &&
\frac 1N \sum_{x\in \bb T_N} (\bb L_N H_\lambda^N) (x/N)\, \eta(x)
\\
\!\!\!\!\!\!\!\!\!\!\!\!\!\! && \quad
+\; \frac a N \sum_{x\in \bb T_N} \big\{ (\bb L_N H_\lambda^N)
(x+1/N) + (\bb L_N H_\lambda^N) (x/N) \big\} \,
(\tau_x h_1) (\eta) \\
\!\!\!\!\!\!\!\!\!\!\!\!\!\! && \qquad
- \; \frac a N \sum_{x\in \bb T_N} (\bb L_N H_\lambda^N)
(x/N) (\tau_x h_2) (\eta)\;,
\end{eqnarray*}
where $\{\tau_x: x\in \bb Z\}$ is the group of translations so that
$(\tau_x \eta)(y) = \eta(x+y)$ for $x$, $y$ in $\bb Z$ and the sum is
understood modulo $N$. Also, $h_1$, $h_2$ are the cylinder functions
\begin{equation*}
h_1(\eta) \;=\; \eta(0) \eta(1)\;,\quad h_2(\eta) \;=\; \eta(-1)
  \eta(1)\; .
\end{equation*}
Since $H_\lambda^N$ is the solution of the resolvent equation
\eqref{f04}, we may replace $\bb L_N H_\lambda^N$ by $U_\lambda^N
= \lambda H_\lambda^N - H$ in the previous formula. In particular,
for all $0\le s<t\le T$,
\begin{equation*}
\Big\vert \int_s^t dr \, N^2 L_N \<\pi^N_r ,H_\lambda^N \> \Big\vert
\;\le\; \frac {(1+3|a|)(t-s)}N \sum_{x\in \bb T_N} |U_\lambda^N (x/N)|
\;.
\end{equation*}
It follows from the first estimate in \eqref{f05} and from Schwarz
inequality that the right hand side is bounded above by $C(H,a) (t-s)$
uniformly in $N$, where $C(H,a)$ is a finite constant depending only
on $a$ and $H$. This proves that the additive part of the
decomposition \eqref{f10} is tight for the uniform topology and
therefore that the sequence of processes $\{\Pi^{\lambda,N}_t :N\ge
1\}$ is tight.

\begin{lemma}
\label{s06}
The sequence of measures $\{\bb Q_{\mu^N}^{W,N} : N\ge 1\}$ is tight
for the uniform topology.
\end{lemma}

\begin{proof}
It is enough to show that for every smooth function $H:\bb T\to\bb R$ 
and every $\epsilon>0$, there exists $\lambda>0$ such that
\begin{equation*}
\lim_{N\to\infty} \bb P_{\mu^N} \Big[
\sup_{0\le t\le T} |\, \Pi^{\lambda,N}_t (\lambda H) - 
\<\pi^N_t, H\>\, | > \epsilon
\Big] \;=\;0
\end{equation*}
because in this case the tightness of $\pi^N_t$ follows from the
tightness of $\Pi^{\lambda,N}_t$.  Since there is at most one particle
per site the expression inside the absolute value is less than or
equal to
\begin{equation*}
\frac{1}{N} \sum_{x \in \bb T_N} \big|\lambda H_\lambda^N (x/N)
- H (x/N)\big|\;.
\end{equation*}
By \eqref{f23} this expression vanishes as $N\uparrow\infty$,
$\lambda\uparrow\infty$.
\end{proof}

\subsection{Uniqueness of limit points} 
\label{ss2}

We prove in this subsection that all limit points $\bb Q^*$ of the
sequence $\bb Q^{W,N}_{\mu_N}$ are concentrated on absolutely
continuous trajectories $\pi(t,du) = \rho(t,u) du$, whose density
$\rho(t,u)$ is a weak solution of the hydrodynamic equation
\eqref{g03} with $l=0$, $r=1$, $\gamma = \rho_0$ and $\Phi(\theta) =
\theta + a\theta^2$.

Let $\bb Q^*$ be a limit point of the sequence $\bb Q^{W,N}_{\mu_N}$
and assume, without loss of generality, that $\bb Q^{W,N}_{\mu_N}$
converges to $\bb Q^*$.

Since there is at most one particle per site, it is clear that $\bb
Q^*$ is concentrated on trajectories $\pi_t(du)$ which are absolutely
continuous with respect to the Lebesgue measure, $\pi_t(du) =
\rho(t,u) du$, and whose density $\rho$ is non-negative and bounded by
$1$.

Fix a function $H: \bb T\to \bb R$ continuously differentiable and
$\lambda>0$.  Recall the definition of the martingale
$M^{N,\lambda}_t$ introduced in the previous section. By \eqref{f02},
for every $\delta>0$,
\begin{equation*}
\lim_{N\to\infty} \bb P_{\mu_N} \Big[ \sup_{0\le t\le T}
\big\vert M^{N,\lambda}_t \big\vert \, > \, \delta \Big] 
\;=\; 0\;.
\end{equation*}

The martingale $M^{N,\lambda}_t$ can be written in terms of the
empirical measure as 
\begin{eqnarray*}
\<\pi^N_t, G^N_\lambda H \> \;-\; \<\pi^N_0, G^N_\lambda H \> \;-\;
\int_0^t ds \, N^2 L_N \<\pi^N_s , G^N_\lambda H \> \;.
\end{eqnarray*}
Therefore, for fixed $0<t\le T$ and $\delta>0$,
\begin{equation*}
\lim_{N\to\infty} \bb Q^{W,N}_{\mu_N} \Big[ \,  
\Big\vert \<\pi^N_t, G^N_\lambda H \> \;-\; 
\<\pi^N_0, G^N_\lambda H \> \;-\;
\int_0^t ds \, N^2 L_N \<\pi^N_s , G^N_\lambda H \>
\Big\vert \, > \, \delta \Big] \;=\; 0\;.
\end{equation*}

Since there is at most one particle per site, by \eqref{f24}, we may
replace $G^N_\lambda H$ by $G_\lambda H$ in the expression $\<\pi^N_t,
G^N_\lambda H \>$, $\<\pi^N_0, G^N_\lambda H \>$ above. 

On the other hand, the expression $N^2 L_N \<\pi^N_s , G^N_\lambda H
\>$ has been computed in the previous subsection. Recall that $\bb L_N
G^N_\lambda H = \lambda G^N_\lambda H - H$. As before, we may replace
$G^N_\lambda H$ by $G_\lambda H$. Let $U_\lambda = \lambda G_\lambda H
- H$.  Since $E_{\nu_\alpha}[h_j] = \alpha^2$, $j=1$, $2$, in view of
\eqref{f05} and by Corollary \ref{s02}, for every $t>0$, $\lambda>0$,
$\delta>0$, $j=1$, $2$,
\begin{equation*}
\lim_{\varepsilon \to 0} \limsup_{N\to\infty} 
\bb P_{\mu_N} \Big[ \, \Big| \int_0^t  ds\, \frac 1N
\sum_{x\in \bb T_N} U_\lambda  (x/N) \Big\{ \tau_x h_j (\eta_s) - 
\big[\eta^{\varepsilon N}_s(x)\big]^2 \Big\} \, \Big|
\, > \, \delta \, \Big]  \;=\; 0\;.
\end{equation*}

Since $\eta^{\varepsilon N}_s(x) = \varepsilon^{-1} \pi^N_s ([x/N, x/N
+ \varepsilon])$, we obtain from the previous considerations that
\begin{eqnarray*}
\!\!\!\!\!\!\!\!\!\!\!\!\!\! &&
\lim_{\varepsilon \to 0} \limsup_{N\to\infty} \bb Q^{W,N}_{\mu_N} \Big[ \,  
\Big\vert \<\pi^N_t, G_\lambda H \> \;-\; \\
\!\!\!\!\!\!\!\!\!\!\!\!\!\! && \qquad\qquad
-\; \<\pi^N_0, G_\lambda H \> \;-\;
\int_0^t ds \, \Big\< \Phi\big (\varepsilon^{-1} \pi^N_s ([\cdot, \cdot
+ \varepsilon]) \big) \,,\, U_\lambda\Big> 
\Big\vert > \delta \Big] \;=\; 0\;.
\end{eqnarray*}

Since $H$ is a smooth function, $G_\lambda H$ and $U_\lambda$ can be
approximated in $L^1(\bb T)$ by continuous functions. Since we assumed
that $\bb Q^{W,N}_{\mu_N}$ converges in the uniform topology to $\bb
Q^*$, we have that
\begin{eqnarray*}
\!\!\!\!\!\!\!\!\!\!\!\!\!\! &&
\lim_{\varepsilon \to 0}  \bb Q^{*} \Big[ \,  
\Big\vert \<\pi_t, G_\lambda H \> \;-\; \\
\!\!\!\!\!\!\!\!\!\!\!\!\!\! && \qquad\qquad\qquad
-\; \<\pi_0, G_\lambda H \> \;-\;
\int_0^t ds \, \Big\< \Phi\big (\varepsilon^{-1} \pi_s ([\cdot, \cdot
+ \varepsilon]) \big) \,,\, U_\lambda\Big> 
\Big\vert > \delta \Big] \;=\; 0\;.
\end{eqnarray*}

As $\bb Q^{*}$ is concentrated on absolutely continuous paths
$\pi_t(du) = \rho(t,u) du$ with positive density bounded by $1$,
$\varepsilon^{-1} \pi_s ([\cdot, \cdot + \varepsilon])$ converges in
$L^1(\bb T)$ to $\rho(s,u)$ as $\varepsilon\downarrow 0$. Thus,
\begin{eqnarray*}
\bb Q^{*} \Big[ \,  
\Big\vert \<\pi_t, G_\lambda H \> \;-\; 
 \<\pi_0, G_\lambda H \> \;-\;
\int_0^t ds \, \< \Phi (\rho_s) \,,\, \mc L_W G_\lambda H \> 
\Big\vert > \delta \Big] \;=\; 0
\end{eqnarray*}
because $U_\lambda = \mc L_W G_\lambda H$. Letting $\delta\downarrow
0$, we see that $\bb Q^{*}$ a.s. 
\begin{eqnarray*}
\<\pi_t, G_\lambda H \> \;-\; \<\pi_0, G_\lambda H \> \;=\;
\int_0^t ds \, \< \Phi (\rho_s) \,,\, \mc L_W G_\lambda H \> \;.
\end{eqnarray*}
This identity can be extended to a countable set of times $t$. Taking
this set to be dense, by continuity of the trajectories $\pi_t$, we
obtain that it holds for all $0\le t\le T$. In the same way, it holds
for any countable family of continuous functions. Taking a countable
set of continuous functions, dense for the uniform topology, we extend
this identity to all continuous function $H$ because $G_\lambda H_n$
converges to $G_\lambda H$ in $L^1(\bb T)$ if $H_n$ converges to $H$
in the uniform topology.  Similarly, we can show that it holds for all
$\lambda>0$, since, for any continuous function $H$, $G_{\lambda_n} H$
converges to $G_\lambda H$ in $L^1(\bb T)$, as $\lambda_n \to
\lambda$.

\begin{proof}[Proof of Proposition \ref{s15}]
In the previous subsection we showed that the sequence of probability
measures $\bb Q^{W,N}_{\mu_N}$ is tight for the uniform topology. We
just proved that all limit points of this sequence are concentrated on
weak solutions of the parabolic equation \eqref{g03}. The statement of
the proposition follows from the uniqueness of weak solutions proved
in Section \ref{sec6}.
\end{proof}

\subsection{Replacement lemma} 

Denote by $H_N (\mu_N | \nu_\alpha)$ the entropy of a probability
measure $\mu_N$ with respext to a stationary state $\nu_\alpha$. We
refer to \cite[Section A1.8]{kl} for a precise definition. By the
explicit formula given in \cite[Theorem A1.8.3]{kl}, we see that there
exists a finite constant $K_0$, depending only on $\alpha$, such that 
\begin{equation}
\label{f06}
H_N (\mu_N | \nu_\alpha) \;\le\; K_0 N
\end{equation}
for all measures $\mu_N$.

Denote by $\< \cdot, \cdot \>_{\nu_\alpha}$ the scalar product of
$L^2(\nu_\alpha)$ and denote by $I^\xi_N$ the convex and lower
semicontinuous \cite[Corollary A1.10.3]{kl} functional defined by
\begin{equation*}
I^\xi_N (f) \;=\; \< - L_N \sqrt f \,,\, \sqrt f\>_{\nu_\alpha}\; ,
\end{equation*}
for all probability densities $f$ with respect to $\nu_\alpha$ (i.e.,
$f\ge 0$ and $\int f d\nu_\alpha =1$). An elementary computation shows
that 
\begin{eqnarray*}
\!\!\!\!\!\!\!\!\!\!\!\!\!\! &&
I^\xi_N (f) \;=\; \sum_{x\in \bb T_N} I^\xi_{x,x+1} (f)\;,
\quad\text{where}\quad \\
\!\!\!\!\!\!\!\!\!\!\!\!\!\! && \qquad
I^\xi_{x,x+1} (f) \;=\; (1/2)\, \xi_x 
\int c_{x,x+1} (\eta) \big\{ \sqrt{f(\sigma^{x,x+1} \eta)} -
\sqrt{f(\eta)} \big\}^2 \, d\nu_\alpha \;.  
\end{eqnarray*}
By \cite[Theorem A1.9.2]{kl}, if $\{S^N_t : t\ge 0\}$ stands for the
semi-group associated to the generator $N^2L_N$, 
\begin{equation*}
H_N (\mu_N S^N_t | \nu_\alpha) \; +\; N^2 \, \int_0^t 
I^\xi_N (f^N_s) \, ds  \;\le\; H_N (\mu_N | \nu_\alpha)\;,
\end{equation*}
provided $f^N_s$ stands for the Radon-Nikodym derivative of $\mu_N
S^N_s$ with respect to $\nu_\alpha$. 
\medskip

For a local function $g: \{0,1\}^{\bb Z} \to \bb R$, let
$\tilde g :[0,1]\to \bb R$ be the expected value of $g$ under the
stationary states:
\begin{equation*}
\tilde g (\alpha) \;=\; E_{\nu_\alpha} [ g(\eta)]\;.
\end{equation*}
For $\ell \ge 1$, let $\eta^\ell (x)$ be the density of particles on
the interval $\{x, \dots, x+ \ell -1\}$:
\begin{equation*}
  \eta^\ell (x) \;=\; \frac 1{\ell}  \sum_{y=x}^{x+\ell-1} \eta(y)\;.
\end{equation*}

\begin{lemma}
\label{s01}
Fix a function $F: N^{-1} \bb T_N \to \bb R$. There exists a finite constant
$C_0$, depending only on $a$, $g$ and $W$, such that
\begin{eqnarray*}
\!\!\!\!\!\!\!\!\!\!\!\! &&
\frac 1N \sum_{x\in \bb T_N} F(x/N) \int  \{ \tau_x g (\eta) - \tilde
g (\eta^{\varepsilon N}(x)) \}\,  f(\eta) \nu_\alpha(d\eta) \\
\!\!\!\!\!\!\!\!\!\!\!\! && \quad
\le\; \frac {C_0} {\varepsilon N^2} \sum_{x\in \bb T_N} \big| F(x/N) \big|
\;+\; \frac {C_0 \varepsilon}{\delta N} \sum_{x\in \bb T_N} F(x/N)^2
 \;+\; \delta N I^\xi_N(f)
\end{eqnarray*}
for all $\delta>0$ and all probability density $f$ with respect to
$\nu_\alpha$. 
\end{lemma}

\begin{proof}
Any local function can be written as a linear combination of functions
of type $\prod_{x\in A} \eta(x)$, for finite sets $A's$. It is
therefore enough to prove the lemma for such functions. We prove the
result for $g(\eta) = \eta(0) \eta(1)$. The general case can be
handled in a similar way.

We estimate first
\begin{equation}
\label{f03}
\frac 1N \sum_{x\in \bb T_N} F(x/N) \int \eta(x) \Big\{ \eta(x+1) -
\frac 1{\varepsilon N}  \sum_{y=x}^{x+\varepsilon N-1} \eta(y) 
\Big\} f(\eta) \nu_\alpha(d\eta) 
\end{equation}
in terms of the functional $I^\xi_N(f)$. The integral can be
rewritten as
\begin{equation*}
\frac 1{\varepsilon N} \sum_{y=x+2}^{x+\varepsilon N-1} 
\sum_{z=x+1}^{y-1} \int \eta(x) 
\{ \eta(z) - \eta(z+1) \} f(\eta) \nu_\alpha(d\eta) \;+\; 
O(\frac 1{\varepsilon N})\;,
\end{equation*}
where the remainder comes from the contribution $y=x$. Writing last
integral as twice the same expression and performing the change of
variables $\eta' = \sigma^{z,z+1}\eta$ in one of them, the previous
integral becomes
\begin{equation*}
(1/2) \int \eta(x) \{ \eta(z) - \eta(z+1) \} \, 
\big\{ f(\eta) - f(\sigma^{z,z+1}\eta) \big\} \, \nu_\alpha(d\eta) \;.
\end{equation*}
Since $a-b = (\sqrt a - \sqrt b)(\sqrt a + \sqrt b)$, by Schwarz
inequality the previous expression is less than or equal to
\begin{eqnarray*}
\!\!\!\!\!\!\!\!\!\!\!\!\!\! &&
\frac A {16 (1-2a^-) \xi_z} \int \eta(x) \{ \eta(z) - \eta(z+1) \}^2 \, 
\big\{ \sqrt{f(\eta)} + \sqrt{f(\sigma^{z,z+1}\eta)} \big\}^2 
\, \nu_\alpha(d\eta) \\
\!\!\!\!\!\!\!\!\!\!\!\!\!\! && \quad \;+\;
\frac {\xi_z} A  \int c_{z,z+1}(\eta) 
\big\{ \sqrt{f(\eta)} - \sqrt{f(\sigma^{z,z+1}\eta)} \big\}^2 
\, \nu_\alpha(d\eta)
\end{eqnarray*}
for every $A>0$. In this formula we used the fact that $c_{z,z+1}$ is
bounded below by $1-2a^-$. Since $f$ is a density with respect to
$\nu_\alpha$, the first expression is bounded by $A/4 (1-2a^-) \xi_z$,
while the second one is equal to $2 A^{-1} I^\xi_{z,z+1}(f)$. Adding
together all previous estimates, we obtain that \eqref{f03} is less
than or equal to
\begin{equation*}
\frac 1{\varepsilon N^2} \sum_{x\in \bb T_N} \big| F(x/N) \big|
\;+\; \frac A{4 (1-2a^-) N} \sum_{x\in \bb T_N} F(x/N)^2
\sum_{z=x+1}^{x+\varepsilon N} \xi_z^{-1} \;+\; \frac {2\varepsilon} A
\sum_{z\in \bb T_N} I^\xi_{z,z+1}(f)\;.
\end{equation*}
By definition of the sequence $\{\xi_z\}$, $\sum_{x+1 \le z\le
  x+\varepsilon N} \xi_z^{-1} \le N [W(1) - W(0)]$. Thus, choosing $A=
2 \varepsilon N^{-1} \delta^{-1}$, for some $\delta>0$, we obtain that
the previous sum is bounded above by
\begin{equation*}
\frac 1{\varepsilon N^2} \sum_{x\in \bb T_N} \big| F(x/N) \big|
\;+\; \frac {C_0 \varepsilon}{\delta N} \sum_{x\in \bb T_N} F(x/N)^2
 \;+\; \delta N I^\xi_N(f)\;.
\end{equation*}

Up to this point we have replaced $\eta(x) \eta(x+1)$ by $\eta(x)
\eta^{\varepsilon N} (x)$. The same arguments permit to replace this
latter expression by $[\eta^{\varepsilon N} (x)]^2$, which concludes
the proof of the lemma.
\end{proof}

\begin{corollary}
\label{s02}
Fix a cylinder function $g$ and a sequence of functions $\{F_N : N\ge
1\}$, $F_N : N^{-1} \bb T_N\to \bb R$ such that
\begin{equation*}
\limsup_{N\to\infty} \frac 1N \sum_{x\in \bb T_N} F_N(x/N)^2 \;<\;
\infty\; . 
\end{equation*}
Then, for any $t>0$ and any sequence of probability measures $\{\mu_N
: N\ge 1\}$ on $\{0,1\}^{\bb T_N}$,
\begin{equation*}
\limsup_{\varepsilon\to 0} \limsup_{N\to\infty}
\bb E_{\mu_N} \Big[ \, \Big| \int_0^t  \frac 1N
\sum_{x\in \bb T_N} F_N(x/N) \, \big \{ \tau_x g (\eta_s) - \tilde g
(\eta^{\varepsilon N}_s(x)) \big \} \Big| \, \Big] \;=\; 0\;.
\end{equation*}
\end{corollary}

\begin{proof}
Fix $0<\alpha <1$. By the entropy and Jensen inequalities, the
expectation appearing in the statement of the lemma is bounded above
by 
\begin{equation*}
\frac {H_N (\mu_N | \nu_\alpha)}{\gamma N} \;+\;
\frac 1{\gamma N} \log \bb E_{\nu_\alpha} \Big[ 
\exp\Big\{ \gamma \Big| \int_0^t  
\sum_{x\in \bb T_N} F_N(x/N) \, \big \{ \tau_x g (\eta_s) - \tilde g
(\eta^{\varepsilon N}_s(x)) \big \} \, \Big| \, \Big\} \, \Big]
\end{equation*}
for all $\gamma >0$. In view of \eqref{f06}, to prove the corollary it
is enough to show that the second term vanishes as $N\uparrow\infty$
and then $\varepsilon\downarrow 0$ for every $\gamma>0$. We may remove
the absolute value inside the exponential because $e^{|x|} \le e^x +
e^{-x}$ and because $\limsup_{N\to\infty} N^{-1} \log\{a_N + b_N\} \le
\max\{ \limsup_{N\to\infty} N^{-1} \log a_N , \limsup_{N\to\infty}
N^{-1}$ $\log b_N \}$. Thus, to prove the corollary, we need to show
that
\begin{equation*}
\limsup_{\varepsilon\to 0} \limsup_{N\to\infty}
\frac 1{N} \log \bb E_{\nu_\alpha} \Big[ 
\exp\Big\{ \gamma \int_0^t  
\sum_{x\in \bb T_N} F_N(x/N) \{ \tau_x g (\eta_s) - \tilde g
(\eta^{\varepsilon N}_s(x)) \} \Big\} \, \Big] \;=\; 0
\end{equation*}
for every $\gamma>0$.

By Feynman-Kac formula, for each fixed $N$ the previous expression is
bounded above by
\begin{equation*}
t\, \sup_{f} \Big\{ \int  \frac \gamma N
\sum_{x\in \bb T_N} F_N(x/N) \{ \tau_x g (\eta) - \tilde g
(\eta^{\varepsilon N}(x)) \} f (\eta) \, d \nu_\alpha\;
-\; N  \hat I^\xi_N (f) \Big\}\;,
\end{equation*}
where  the supremum  is carried  over all  density functions  $f$ with
respect to  $\nu_\alpha$. Letting $\delta  =1$ in Lemma  \ref{s01}, we
obtain that the previous expression is less than or equal to
\begin{equation*}
\frac {C_0 \gamma} {\varepsilon N^2} 
\sum_{x\in \bb T_N} \big| F_N(x/N) \big|
\;+\; \frac {C_0 \gamma \varepsilon}{N} 
\sum_{x\in \bb T_N} F_N(x/N)^2
\end{equation*}
for some finite constant $C_0$ which depends on $g$ and $W$. By
assumption on the sequence $\{F_N\}$, for every $\gamma>0$, this
expression vanishes as $N\uparrow\infty$ and then
$\varepsilon\downarrow 0$. This concludes the proof of the lemma.
\end{proof}

\section{Energy estimate}
\label{sec4}

We prove in this section that any limit point $\bb Q^*_{W}$ of the
sequence $\bb Q_{\mu_N}^{W,N}$ is concentrated on trajectories
$\rho(t,u) du$ with finite energy. Though not needed in the proof of
the law of large numbers of the empirical measure $\pi^N$, this
estimate plays an important role in the proof of the large deviations
principle. 

Let $\bb Q^*_{W}$ be a limit point of the sequence $\bb
Q_{\mu_N}^{W,N}$ and assume without loss of generality that the
sequence $\bb Q_{\mu_N}^{W,N}$ converges to $\bb Q^*_{W}$.  Denote by
$\partial_u$ the partial derivative of a function with respect to the
space variable. Let $L^2_W([0,T]\times \bb T)$ be the Hilbert space of
measurable functions $H: [0,T]\times \bb T\to\bb R$ such that
\begin{equation*}
\int_0^Tds \int_{\bb T} dW(u) \, H (s, u)^2 \;<\; \infty\;,
\end{equation*}
endowed with the scalar product $\<\!\< H,G \>\!\>_W$ defined by
\begin{equation*}
\<\!\< H,G \>\!\>_W \;=\; \int_0^Tds \int_{\bb T} dW(u) \, H (s, u)
\, G(s,u)\;.
\end{equation*}

\begin{proposition}
\label{s05}
The measure $\bb Q^*_{W}$ is concentrated on paths $\rho(t,u) du$ with
the property that there exists a function in $L^2_W([0,T]\times \bb
T)$, denoted by $d\Phi/dW$, such that
\begin{equation*}
\int_0^Tds \int_{\bb T} du \, (\partial_u H) (s, u) \, \Phi(\rho(s,u))
\;=\; -\; \int_0^Tds \int_{\bb T} dW(u) \, (d\Phi/dW) (s, u) \, H (s, u)
\end{equation*}
for all functions $H$ in $C^{0,1}([0,T]\times \bb T)$.
\end{proposition}

The previous result follows from the next lemma.  Recall the
definition of the constant $K_0$ given in \eqref{f06}.

\begin{lemma}
\label{s03}
There exists a finite constant $K_1$, depending only on $a$, such that
\begin{eqnarray*}
\!\!\!\!\!\!\!\!\!\!\!\!\! &&
E_{\bb Q^*_{W}} \Big[ \sup_H \Big\{ \int_0^T ds\, \int_{\bb T}
du \, (\partial_u H) (s, u) \, \Phi(\rho(s,u)) \\
\!\!\!\!\!\!\!\!\!\!\!\!\! && \qquad\qquad\qquad\qquad\qquad
\qquad\qquad\qquad
- \; K_1 \int_0^T ds\, \int_{\bb T} H (s, u)^2 
\, dW(u) \Big\} \Big] \; \le \; K_0 \; ,
\end{eqnarray*}
where the supremum is carried over all functions $H$ in
$C^{0,1}([0,T]\times \bb T)$.
\end{lemma}

\begin{proof}[Proof of Proposition \ref{s05}]
Denote by $\ell : C^{0,1}([0,T]\times \bb T) \to \bb R$ the linear
functional defined by
\begin{equation*}
\ell (H) \;=\; \int_0^T ds\, \int_{\bb T}
du \, (\partial_u H) (s, u) \, \Phi(\rho(s,u))\;.
\end{equation*}
Since $C^{0,1}([0,T]\times \bb T)$ is dense in $L^2_W([0,T]\times \bb
T)$, by Lemma \ref{s03}, $\ell$ is $\bb Q^*_{W}$-almost surely finite
in $L^2_W([0,T]\times \bb T)$. In particular, by Riesz representation
theorem, there exists a function $G$ in $L^2_W([0,T]\times \bb T)$
such that
\begin{equation*}
\ell (H) \;=\; - \int_0^T ds\, \int_{\bb T}
dW(u) \, H (s, u) \, G(s,u)\;.
\end{equation*}
This concludes the proof of the proposition.
\end{proof}

The proof of Lemma \ref{s03} relies on the following result.  For a
smooth function $H\colon \bb T\to \bb R$, $\delta >0$, $\varepsilon
>0$ and a positive integer $N$, define $W_N(\varepsilon, \delta, H,
\eta)$ by
\begin{eqnarray*}
W_N(\varepsilon, \delta, H, \eta ) &=&
\sum_{x\in\bb T_N} H(x/N) \frac 1{\varepsilon N} 
\, \Big\{ \Phi(\eta^{\delta N} (x)) -
\Phi(\eta^{\delta N} (x + \varepsilon N)) \Big\} \\
&-& \frac {K_1}{\varepsilon N} \sum_{x\in\bb T_N} 
H(x/N)^2 \{ W([x+\varepsilon N +1]/N) - W(x/N) \}\; .
\end{eqnarray*}

\begin{lemma}
\label{s04}
Consider a sequence $\{H_\ell,\, \ell\ge 1\}$ dense in
$C^{0,1}([0,T]\times \bb T)$.  For every $k\ge 1$, and every
$\varepsilon >0$,
\begin{equation*}
\limsup_{\delta\to 0} \limsup_{N\to\infty}
\bb E_{\mu^N} \Big[ \max_{1\le i\le k} \Big\{
\int_0^T W_N(\varepsilon, \delta, H_i (s, \cdot) , \eta_s ) \, 
ds \Big\} \Big] \;\le\; K_0\; .
\end{equation*}
\end{lemma}

\begin{proof}
It follows from the replacement lemma that in order to prove
the lemma we just need to show that
\begin{equation*}
\limsup_{N\to\infty} \bb E_{\mu^N} \Big[ \max_{1\le i\le k} \Big\{
\int_0^T W_N(\varepsilon,  H_i (s, \cdot) , \eta_s ) \, ds \Big\}
\Big] \;\le\; K_0\; ,
\end{equation*}
where
\begin{eqnarray*}
W_N(\varepsilon, H , \eta ) &=&
\frac 1{\varepsilon N} \sum_{x\in\bb T_N} H(x/N) 
\big\{ \tau_x g(\eta) -  \tau_{x + \varepsilon N} g(\eta)\big\} \\
&-&  \frac {K_1}{\varepsilon N} \sum_{x\in\bb T_N} H(x/N)^2 
\{ W([x+\varepsilon N +1]/N) - W(x/N) \}\; ,
\end{eqnarray*}
and $g(\eta) = \eta(0) + 2a \eta(0)\eta(1)$.

By the entropy and the Jensen inequality,
for each fixed $N$, the previous expectation is bounded above by
\begin{equation*}
\frac {H(\mu^N \vert \nu_{\alpha})}{ N} \; +\; \frac 1{N}
\log \bb E_{\nu_{\alpha}} \Big[ \exp\Big\{
\max_{1\le i\le k} \Big\{ N \int_0^T ds\,
W_N(\varepsilon,  H_i (s, \cdot) , \eta_s ) \Big\} \Big\} \Big] \; .
\end{equation*}
By \eqref{f06}, the first term is bounded by $K_0$.  Since $\exp\{
\max_{1\le j\le k} a_j \}$ is bounded above by $\sum_{1\le j\le k}
\exp\{a_j\}$ and since $\limsup_N N^{-1} \log \{a_N + b_N\}$ is less
than or equal to the maximum of $\limsup_N N^{-1} \log a_N$ and
$\limsup_N N^{-1} \log b_N$, the limit, as $N\uparrow\infty$, of the
second term of the previous expression is less than or equal to
\begin{equation*}
\max_{1\le i \le k} \limsup_{N\to\infty} \frac 1{N} \log
\bb E_{\nu_{\alpha}} \Big[ \exp
\Big\{ N  \int_0^T ds\,  W_N(\varepsilon,  H_i (s, \cdot) , \eta_s )
\Big\} \Big] \; .
\end{equation*}
We now prove that for each fixed $i$ the above limit is nonpositive.

Fix $1\le i\le k$.
By Feynman--Kac formula and the variational formula for the
largest eigenvalue of a symmetric operator, for each fixed
$N$, the previous expression is bounded above by
\begin{equation*}
\int_0^T ds\, \sup_{f} \Big\{  \int W_N(\varepsilon,
H_i (s, \cdot) , \eta )
f(\eta) \nu_{\alpha} (d\eta) - N I^\xi_N (f) \Big\}\; .
\end{equation*}
In this formula the supremum is taken over all probability densities
$f$ with respect to $\nu_{\alpha}$.

It remains to rewrite $\eta(x) \eta(x+1) - \eta(x+\varepsilon
N)\eta(x+ \varepsilon N +1)$ as $\eta(x) \{\eta(x+1) - \eta(x+
\varepsilon N +1)\} + \eta(x+ \varepsilon N +1) \{\eta(x) -
\eta(x+\varepsilon N)\}$ and to repeat the arguments presented in the
proof of Lemma \ref{s01} to conclude.
\end{proof}

\begin{proof}[Proof of Lemma \ref{s03}] 
Assume without loss of generality that $\bb Q_{\mu_N}^{W,N}$ converges
to $\bb Q^*_{W}$. Consider a sequence $\{H_\ell,\, \ell\ge 1\}$ dense in
$C^{0,1}([0,T]\times \bb T)$. By Lemma \ref{s04}, for every $k\ge 1$
\begin{eqnarray*}
\!\!\!\!\!\!\!\!\!\!\!\!\! &&
\limsup_{\delta\to 0} E_{\bb Q^*_{W}}\Big[ \max_{1\le i\le k} 
\Big\{ \frac 1{\varepsilon} \int_0^T ds\, \int_{\bb T} du \,
H_i (s,u) \, \Big\{ \Phi(\rho^\delta_s (u)) -
\Phi(\rho^\delta_s (u + \varepsilon)) \Big\} \\
\!\!\!\!\!\!\!\!\!\!\!\!\! &&
\qquad\qquad -\; \frac {K_1} {\varepsilon} \int_0^T ds\, 
\int_{\bb T} du \, H_i(s,u)^2  \, [W(u+\varepsilon) - W(u)]
\Big\} \Big]\; \le \; K_0\; ,
\end{eqnarray*}
where $\rho^\delta_s (u) = (\rho_s * \iota_\delta)(u)$ and
$\iota_\delta$ is the approximation of the identity $\iota_\delta
(\cdot) = (2\delta)^{-1} \mb 1\{ [-\delta , \delta]\} (\cdot)$.

Letting $\delta\downarrow 0$, changing variables and then letting
$\varepsilon\downarrow 0$, we obtain that
\begin{eqnarray*}
\!\!\!\!\!\!\!\!\!\!\!\!\! &&
E_{\bb Q^*_{W}}\Big[ \max_{1\le i\le k} \Big\{
\int_0^T ds\, \int_{\bb T} (\partial_u H_i) (s,u)
\Phi(\rho (s,u)) \, du \\
\!\!\!\!\!\!\!\!\!\!\!\!\! && \qquad\qquad\qquad\qquad\qquad\qquad
- \;K_1 \int_0^T ds\, \int_{\bb T} H_i(s,u)^2 dW(u)
\Big\} \Big] \;\le \; K_0\; .
\end{eqnarray*}
To conclude the proof it remains to apply the monotone convergence
theorem and recall that $\{H_\ell, \, \ell\ge 1\}$ is a dense sequence
in $C^{0,1}([0, T]\times \bb T)$ for the norm $\Vert H\Vert_\infty
+ \Vert (\partial_u H)\Vert_\infty$.
\end{proof}

\section{Uniqueness of weak solutions of \eqref{g03}}
\label{sec6}

Recall that we denote by $\<\cdot, \cdot\>$ the inner product of the
Hilbert space $L^2(\bb T)$ and that $\{G_\lambda : \lambda >0\}$
stands for the resolvents associated to $\mc L_W$.

Let $\rho$ be a weak solution of the hydrodynamic equation
\eqref{g03}. Since $\rho$, $\Phi(\rho)$ are bounded, since the smooth
functions are dense in $L^2(\bb T)$ and since $\mc L_W G_\lambda = -
\bb I + \lambda G_\lambda$ are bounded operatores, for any function
$H$ in $L^2(\bb T)$,
\begin{equation*}
\< \rho_t, G_\lambda H\> \;-\; \< \gamma , G_\lambda H\> 
\;=\; \int_0^t \< \Phi(\rho_s) , \mc L_W G_\lambda H \>\, ds
\end{equation*}
for all $t>0$ and all $\lambda>0$.

Let $\rho : \bb R_+ \times \bb T \to [l,r]$ be a weak solution of
\eqref{g03}. We claim that 
\begin{equation}
\label{f08}
\< \rho_t \,,\, G_\lambda \rho_t\> \;-\; \< \rho_0 \,,\, G_\lambda
\rho_0\> \;=\; 2 \int_{0}^{t} \< \Phi(\rho_{s}) \,,\, \mc L_W G_\lambda
\rho_{s}\>\, ds
\end{equation}
for all $t>0$ and $\lambda>0$.

To prove this claim, fix $\lambda>0$, $t>0$ and consider a partition
$0=t_0< t_1 < \cdots < t_n = t$ of the interval $[0,t]$ so that
\begin{eqnarray*}
\< \rho_t \,,\, G_\lambda \rho_t\> \;-\; \< \rho_0 \,,\, G_\lambda
\rho_0\> & =& \sum_{k=0}^{n-1} \< \rho_{t_{k+1}} \,,\, G_\lambda
\rho_{t_{k+1}}\> \;-\; \< \rho_{t_{k+1}} \,,\, G_\lambda
\rho_{t_{k}}\> \\
&+& \sum_{k=0}^{n-1} \< \rho_{t_{k+1}} \,,\, G_\lambda \rho_{t_{k}}\>
\;-\; \< \rho_{t_{k}} \,,\, G_\lambda \rho_{t_{k}}\>\;.
\end{eqnarray*}

We handle the first term, the second one being similar. Since
$G_\lambda$ is self-adjoint in $L^2(\bb T)$, since $\rho_{t_{k+1}}$
belongs to $L^2(\bb T)$ and since $\rho$ is a weak solution of
\eqref{g03},
\begin{equation*}
\< \rho_{t_{k+1}} \,,\, G_\lambda \rho_{t_{k+1}}\> \;-\; 
\< \rho_{t_{k+1}} \,,\, G_\lambda \rho_{t_{k}}\> \;=\;
\int_{t_{k}}^{t_{k+1}} \< \Phi(\rho_{s}) \,,\, \mc L_W G_\lambda
\rho_{t_{k+1}}\>\, ds\;.
\end{equation*}
Add and subtract on the right hand side $\< \Phi(\rho_{s}) \,,\, \mc
L_W G_\lambda \rho_{s}\>$. The time integral of this term is exactly
the expression announced in \eqref{f08} and the remainder is given by
\begin{equation*}
\int_{t_{k}}^{t_{k+1}} \Big\{ \< \Phi(\rho_{s}) \,,\, \mc L_W G_\lambda
\rho_{t_{k+1}}\> \;-\; \< \Phi(\rho_{s}) \,,\, \mc L_W G_\lambda
\rho_{s}\> \Big\} \, ds\;.
\end{equation*}

Since $\mc L_W G_\lambda = - \bb I + \lambda G_\lambda$, where $\bb I$
is the indentity, and since $G_\lambda$ is self-adjoint, we may
rewrite the previous difference as
\begin{equation*}
- \, \Big\{ \< \Phi(\rho_{s}) \,,\, \rho_{t_{k+1}}\> \;-\; 
\< \Phi(\rho_{s}) \,,\, \rho_{s}\>  \Big\} \; 
+\; \lambda \Big\{ \< G_\lambda \Phi(\rho_{s}) \,,\,  
\rho_{t_{k+1}}\> \;-\; \< G_\lambda \Phi(\rho_{s}) \,,\, \rho_{s}\>
\Big\}\;.
\end{equation*}
The time integral between $t_k$ and $t_{k+1}$ of the second term is
equal to
\begin{equation*}
\lambda \int_{t_{k}}^{t_{k+1}} ds \int_{s}^{t_{k+1}}  
\< \mc L_W G_\lambda \Phi(\rho_{s}) \,,\, \Phi(\rho_{r})\> \; dr
\end{equation*}
because $\rho$ is a weak solution of \eqref{g03} and $\Phi(\rho_{s})$
belongs to $L^2(\bb T)$. It follows from the boundness of the operator
$\mc L_W G_\lambda$ and from the boundness of $\Phi (\rho)$ that this
expression is of order $(t_{k+1} - t_{k})^2$.

To conclude the proof of claim \eqref{f08} it remains to show that
\begin{equation*}
\sum_{k=0}^{n-1} \int_{t_{k}}^{t_{k+1}} \Big\{ \< \Phi(\rho_{s}) \,
,\, \rho_{t_{k+1}}\> \;-\;  \< \Phi(\rho_{s}) \,,\, \rho_{s}\>  \Big\}
\, ds
\end{equation*}
vanishes as the mesh of the partition tends to $0$. Fix $\varepsilon
>0$ and choose $\beta$ large enough for
\begin{equation*}
\int_0^t ds \int_{\bb T}  \Big\{ \beta G_\beta \Phi(\rho(s,u)) 
- \Phi (\rho(s,u)) \Big\}^2  \,du \;\le\; \varepsilon\;.
\end{equation*}
This is possible because $\Phi(\rho)$ is bounded, $\{\beta G_\beta :
\beta>0\}$ are uniformly bounded operators, and $\beta G_\beta
\Phi(\rho(s,\cdot))$ converges to $\Phi(\rho(s,\cdot))$ in $L^2(\bb
T)$, as $\beta\uparrow\infty$, for all $0\le s\le t$.

Paying a price of order $\sqrt{\varepsilon}$, because $l\le\rho\le r$,
we may replace $\Phi(\rho_{s})$ in the penultimate formula by $\beta
G_\beta \Phi(\rho_{s})$. After this replacement, since $\rho$ is weak
solution, we may rewrite the sum as
\begin{equation*}
\beta \sum_{k=0}^{n-1} \int_{t_{k}}^{t_{k+1}} ds \int_s^{t_{k+1}}
\< \mc L_W G_\beta \Phi(\rho_{s}) \,,\, \Phi(\rho_{r}) \> 
\; dr\;.
\end{equation*}
We have already seen that this expression vanishes as the mesh of the
partition tends to $0$. This proves \eqref{f08}.

Recall the definition of the constant $B$ given at the beginning of
Subsection \ref{ss2.3}

\begin{lemma}
\label{s12}
Fix two density profiles $\gamma^1$, $\gamma^2:\bb T\to [l,r]$ and
denote by $\rho^1$, $\rho^2$ weak solutions of \eqref{g03} with
initial value $\gamma^1$, $\gamma^2$, respectively. Then, 
\begin{equation*}
\Big \< \rho^1_t - \rho^2_t \,,\, G_\lambda 
\big[ \rho^1_t - \rho^2_t\big] \Big \> \;\le\; 
\Big \< \gamma^1 - \gamma^2 \,,\, G_\lambda 
\big[ \gamma^1 - \gamma^2\big] \Big \> \, e^{ B \lambda t/2}
\end{equation*}
for all $\lambda>0$, $t>0$.  In particular, there exists at most one
weak solution of \eqref{g03}.
\end{lemma}

\begin{proof}
Fix two density profiles $\gamma^1$, $\gamma^2:\bb T\to [l,r]$.  Let
$\rho^1$, $\rho^2$ be two weak solutions with initial value
$\gamma^1$, $\gamma^2$, respectively. By \eqref{f08}, for any
$\lambda>0$,
\begin{eqnarray*}
\!\!\!\!\!\!\!\!\!\!\!\!\!\!\! &&
\Big \< \rho^1_t - \rho^2_t \,,\, G_\lambda 
\big[ \rho^1_t - \rho^2_t\big] \Big \> \;
-\; \Big \< \gamma^1 - \gamma^2 \,,\, G_\lambda 
\big[ \gamma^1 - \gamma^2\big] \Big \> \;= \\
\!\!\!\!\!\!\!\!\!\!\!\!\!\!\! && 
- 2 \int_{0}^{t} \< \Phi(\rho^1_{s}) - \Phi(\rho^2_{s}) \,,\, 
\rho^1_{s} -  \rho^2_{s} \>\, ds 
\;+\; 2 \lambda \int_{0}^{t} \Big \< \Phi(\rho^1_{s}) - \Phi(\rho^2_{s}) 
\,,\,  G_\lambda \big[ \rho^1_{s} -  \rho^2_{s} \big]\Big \>\, ds\;.
\end{eqnarray*}

By Schwarz inequality, the second term on the right hand side is
bounded above by
\begin{equation*}
\frac 1A  \int_{0}^{t} \Big \< \Phi(\rho^1_{s}) - \Phi(\rho^2_{s}) 
\,,\,  G_\lambda \big[  \Phi(\rho^1_{s}) - \Phi(\rho^2_{s}) \big]\Big \>\, ds
\;+\; A\lambda^2  \int_{0}^{t} \Big \< \rho^1_{s} -  \rho^2_{s}
\,,\,  G_\lambda \big[ \rho^1_{s} -  \rho^2_{s} \big]\Big \>\, ds
\end{equation*}
for every $A>0$. Since the operator $G_\lambda$ is bounded by
$\lambda^{-1}$, and since $\Phi'$ is bounded by $B$, the first term of
the previous expression is less than or equal to
\begin{equation*}
\frac B{A \lambda} \int_{0}^{t} \Big \< \rho^1_{s} - \rho^2_{s} 
\,,\,    \Phi(\rho^1_{s}) - \Phi(\rho^2_{s}) \Big \>\, ds \;.
\end{equation*}
Choosing $A = B/2\lambda$, this expression cancels with the first term
on the right hand side of the first formula. In particular, the left
hand side of this formula is bounded by
\begin{equation*}
\frac {B \lambda}2 \int_{0}^{t} \Big \< \rho^1_{s} -  \rho^2_{s}
\,,\,  G_\lambda \big[ \rho^1_{s} -  \rho^2_{s} \big]\Big \>\, ds\;.
\end{equation*}
It remains to recall Gronwall's inequality to conclude.
\end{proof}


\begin{thebibliography}{99}
  
\bibitem{dyn2} E.\ B.\ Dynkin, {\em Markov processes}. Volume II.
  Grundlehren der Mathematischen Wissenschaften [Fundamental
  Principles of Mathematical Sciences], 122.  Springer-Verlag, Berlin,
  1965.

\bibitem{fag} A.\ Faggionato {\em Bulk diffusion of 1d exclusion
    process with bond disorder}.  Markov Process. Relat. Fields
  {\bf 13}, 519--542, (2007).


\bibitem{fjl} A. Faggionato, M.\ Jara, C.\ Landim, {\em Hydrodynamic
    behavior of one dimensional subdiffusive exclusion processes with
    random conductances}. arXiv:0709.0306 . To appear in Probab. Th.
  Rel. Fields (2008).
  
\bibitem{f} U. Freiberg {\em Analytical properties of measure
    geometric Krein-Feller-operators on the real line} Math.  Nachr.
  {\bf 260} 34 -- 47, (2003).
  
\bibitem{jl2} M.\ Jara, C.\ Landim, {\em Quenched nonequilibrium
    central limit theorem for a tagged particle in the exclusion
    process with bond disorder}. arXiv:math/0603653.  Ann.  Inst. H.
  Poincar\'e, Probab. Stat. {\bf 44}, 341--361, (2008).

\bibitem{kl} C.\ Kipnis, C.\ Landim, {\em Scaling limits of interacting
  particle systems}. Grundlehren der Mathematischen Wissenschaften
  [Fundamental Principles of Mathematical Sciences], 320.
  Springer-Verlag, Berlin, 1999.

\bibitem{lo1} J.-U.\ L\"obus, {\em  Generalized second order differential
  operators}.  Math. Nachr. {\bf 152}, 229-245 (1991).

\bibitem{lo2} J.-U.\ L\"obus, {\em  Construction and generators of
  one-dimensional quasi-diffusions with applications to selfaffine
  diffusions and brownian motion on the Cantor set}.  Stoch. and
  Stoch. Rep. {\bf 42}, 93--114, (1993).
  
\bibitem{m} P. Mandl, {\em Analytical treatment of one-dimensional
    {M}arkov processes}, Grundlehren der mathematischen
  Wissenschaften, 151.  Springer-Verlag, Berlin, 1968.

\bibitem{n} K.\ Nagy, {\em Symmetric random walk in random
    environment}. Period. Math. Ung. {\bf 45}, 101--120 (2002).

\bibitem{z} E. Zeidler, {\em Applied Functional Analysis. Applications
    to Mathematical Physics.}. Applied Mathematical Sciences, 108.
  Springer-Verlag, New York, 1995.
  
\end{thebibliography}
\end{document}